\documentclass[11pt,reqno, twoside]{amsart}
\usepackage{amscd}
\usepackage{amsfonts}
\usepackage{amsmath}
\usepackage{amssymb}
\usepackage{amsthm}
\usepackage{fancyhdr}
\usepackage{latexsym}
\usepackage[colorlinks=true, pdfstartview=FitV, linkcolor=blue, citecolor=blue, urlcolor=blue]{hyperref}
\usepackage{enumitem}      
\usepackage{mathtools}            
\usepackage{upgreek}
\usepackage{indentfirst} 
\usepackage{schemata} 
\usepackage{blindtext, rotating}   
\usepackage{soul} 
\setstcolor{red}
\usepackage{color}
\usepackage{tikz}
\usetikzlibrary{arrows.meta}
\usetikzlibrary{decorations.markings}
\tikzset{->-/.style={decoration={
  markings,
  mark=at position #1 with {\arrow{>}}},postaction={decorate}}}
  \tikzset{middlearrow/.style={
        decoration={markings,
            mark= at position 0.55 with {\arrow{#1}} ,
        },
        postaction={decorate}
    }
}
\usepackage{caption}           
\newcommand{\eee}[1]{\begin{equation}#1\end{equation}}
\newcommand{\sss}[1]{\begin{subequations}#1\end{subequations}}
\newcommand{\ddd}[1]{\begin{alignat}{2}#1\end{alignat}}

\newcommand{\nn}{\nonumber}
\newcommand{\p}{\partial}

\definecolor{ddgreen}{RGB}{0,170,0}

\renewcommand{\b}{\mathcolor{blue}}

\newcommand{\no}[1]{\left\| #1 \right\|}

\newcommand{\what}{\widehat}
\usepackage{scalerel,stackengine}
\stackMath
\newcommand\wwhat[1]{%
\savestack{\tmpbox}{\stretchto{%
  \scaleto{%
    \scalerel*[\widthof{\ensuremath{#1}}]{\kern-.6pt\bigwedge\kern-.6pt}%
    {\rule[-\textheight/2]{1ex}{\textheight}}
  }{\textheight}%
}{0.5ex}}%
\stackon[1pt]{#1}{\tmpbox}%
}
\makeatletter
\renewcommand\subsection{\@startsection{subsection}{2}%
  \z@{-0.8\linespacing\@plus-0.7\linespacing}{0.7\linespacing}%
  {\normalfont\bfseries}}
\renewcommand\subsubsection{\@startsection{subsubsection}{3}%
  \z@{-0.8\linespacing\@plus-0.7\linespacing}{0.7\linespacing}%
  {\normalfont\itshape}}
\makeatother
\makeatletter
\def\mathcolor#1#{\@mathcolor{#1}}
\def\@mathcolor#1#2#3{%
\protect\leavevmode
\begingroup
\color#1{#2}#3%
\endgroup
}
\makeatother
\theoremstyle{plain}  
\newtheorem{theorem}{Theorem}[section]

\newtheorem{lemma}{Lemma}[section]

\theoremstyle{definition}

\newtheorem{remark}{Remark}[section]
\newenvironment{Proof}[1][\proofname]
{\proof[\textnormal{\textbf{#1.}}]}{\endproof}
\newcommand{\bp}{\begin{Proof}}
\newcommand{\ep}{\end{Proof}}

\DeclareMathSizes{12}{12}{7}{5}
\usepackage[framemethod=TikZ]{mdframed}
\usepackage[breakable, theorems, skins]{tcolorbox}
\tcbset{enhanced}

\numberwithin{figure}{section}
\numberwithin{equation}{section}
\usepackage{geometry}
\makeatletter
\def\l@section{\@tocline{1}{0pt}{1pc}{}{}}
\def\l@subsection{\@tocline{2}{0pt}{1pc}{4.6em}{}}
\def\l@subsubsection{\@tocline{3}{0pt}{1pc}{7.6em}{}}
\renewcommand{\tocsection}[3]{%
  \indentlabel{\@ifnotempty{#2}{\makebox[2.3em][l]{%
    \ignorespaces#1 #2.\hfill}}}#3}
\renewcommand{\tocsubsection}[3]{%
  \indentlabel{\@ifnotempty{#2}{\hspace*{2.3em}\makebox[2.3em][l]{%
    \ignorespaces#1 #2.\hfill}}}#3}
\renewcommand{\tocsubsubsection}[3]{%
  \indentlabel{\@ifnotempty{#2}{\hspace*{4.6em}\makebox[3em][l]{%
    \ignorespaces#1 #2.\hfill}}}#3}
\makeatother 
\begin{document}
\title{Well-Posedness of the Nonlinear Schr\"odinger Equation 
\\
on the Half-Plane
}
\author{A. Alexandrou Himonas$^*$ \& Dionyssios Mantzavinos}
\begin{abstract}
The initial-boundary value problem (IBVP)  for the  nonlinear Schr\"odinger (NLS) equation on the half-plane with nonzero boundary data is studied by advancing a novel approach recently developed for the well-posedness of the cubic NLS on the half-line, which takes advantage of the solution formula produced by the unified transform of Fokas for the associated linear IBVP. 
For initial data in  Sobolev spaces on the half-plane and boundary data in   Bourgain spaces arising naturally when the linear IBVP is solved with zero initial data, the present work provides a local well-posedness result for NLS initial-boundary value problems in higher dimensions.
\end{abstract}
\date{October 4, 2018. \mbox{}$^*$\!\textit{Corresponding author}: himonas.1@nd.edu}
\keywords{
2D nonlinear Schr\"odinger equation, 
initial-boundary value problem,
well-posedness, 
Sobolev spaces,
Bourgain spaces, 
unified transform, Fokas method,   
linear space-time estimates, 
Laplace transform.}

\subjclass[2010]{Primary: 35Q55, 35G16, 35G31}

\maketitle
\markboth
{Well-Posedness of the Nonlinear Schr\"odinger Equation on the Half-Plane}
{A. Himonas \& D. Mantzavinos}
%
%
%
%
%
%
%
%
\section{Introduction and Results}
\label{2d-nls-intro}
We consider the following  initial-boundary value problem (IBVP) for the 
two-dimensional nonlinear Schr\"odinger (NLS) equation  on the half-plane:  
\sss{\label{2d-nls-nls-ibvp}
\ddd{
&i u_t + u_{x_1x_1}+u_{x_2x_2} \pm |u|^{p-1} u=0,
\quad
&& (x_1, x_2, t)\in \mathbb R\times \mathbb R^+ \times (0, T), 
\\
&u(x_1, x_2, 0) = u_0(x_1, x_2), \quad
&& (x_1, x_2)\in \mathbb R\times \overline{\mathbb R^+},
\\
&u(x_1, 0, t) = g_0(x_1, t), && (x_1, t)\in \mathbb R \times [0, T],
}
}
where $p=3, 5, 7, \ldots$ and $T<1$, and prove its well-posedness  for initial data   $u_0(x_1, x_2)$ in the Sobolev space $H^s(\mathbb R_{x_1}\times \mathbb R_{x_2}^+)$, $1<s\leqslant \frac 32$,  and boundary data $g_0(x_1, t)$ in a Bourgain space that arises naturally when the associated linear IBVP is solved with a zero initial datum. 
This result is proven by advancing into two dimensions a novel approach recently introduced for the well-posedness of IBVPs on the half-line \cite{fhm2017, fhm2016, hm2015}.

Remarkably, the Bourgain space associated with the one-dimensional NLS initial value problem (IVP) on the line emerges spontaneously in our work as the natural space for the boundary data. Furthermore,  the regularity in the boundary variables $x_1, t$ of the solution of the linear Schr\"odinger IVP on the plane is described by the aforementioned Bourgain  space,  thereby extending into higher dimensions the result of \cite{kpv1991} on the time regularity of linear Schr\"odinger IVP on the line. 
In this respect, the Bourgain space from the one-dimensional NLS on the line comes forth as the optimal choice for the boundary data space of the two-dimensional NLS on the half-plane.

The NLS equation $iu_t + \Delta u \pm |u|^{p-1}u=0$ is a universal model, arising in multiple areas of mathematical physics such as nonlinear optics \cite{t1964, a2007}, water waves \cite{z1968, ho1972, p1983, css1992},  plasmas \cite{ww1977}, and Bose-Einstein condensates \cite{ps2003, cp2011, kss2011, ch2013}. As such, it has been studied extensively, from many points of view, and in a variety of different contexts. In one dimension, Zakharov and Shabat showed that the cubic case $p=3$ is a completely integrable system \cite{zs1972} (see also \cite{akns1974}). Thus, under the assumption of sufficient smoothness and decay, they were able to study the associated IVP on the line by analyzing the Lax pair of the equation via the inverse scattering transform. For rough data, the NLS IVP has been studied in great detail via harmonic analysis  techniques. In particular, Tsutsumi \cite{t1987} proved global well-posedness in $L^2(\mathbb R^n)$ in the subcritical case $1<p<1+\frac 4n$. Cazenave and Weissler \cite{cw1989} extended this result to the critical case $p = 1 + \frac 4n$ and later to all $p>1$ for data in $H^s(\mathbb R^n)$, $s>0$, using Besov spaces \cite{cw1990}. Other relevant works include Ginibre and Velo \cite{gv1979, gv1985}, Kato \cite{k1987, k1995}, and Constantin and Saut \cite{cs1989}. The periodic case turned out to be more challenging and required the 1993 breakthrough of Bourgain \cite{b1993}, who proved well-posedness for data in $H^s(\mathbb T)$, $s\geqslant 0$, by introducing the celebrated $X^{s, b}$ spaces (with $b=\frac 12$). Further results on well-posedness and ill-posedness in one and higher dimensions for the periodic and the non-periodic IVP are available in \cite{b1998, b1999b, b1999, kpv2001, ckstt2008, kv2012, d2012} and in numerous other works in the literature. 

Contrary to the IVP, progress towards the rigorous analysis of IBVPs for NLS and other nonlinear evolution equations is rather limited. This can be largely attributed to the absence of the Fourier transform and underlying theory in the case of domains with a boundary. Recall, in particular, that the procedure for establishing local well-posedness of a nonlinear problem via a contraction mapping argument is initiated by defining an iteration map through the solution of the associated forced linear problem. In the case of the IVP, this task is straightforward thanks to the Fourier transform and Duhamel's principle. For IBVPs, however, a proper spatial transform is not available and hence a significant challenge is already present at the first step of the analysis. 
This difficulty is reflected in two independent approaches that were introduced in the early 2000s for showing well-posedness of the Korteweg-de Vries (KdV) equation on the half-line, namely the works of Colliander and Kenig \cite{ck2002} (in fact, this work is concerned with the generalized KdV equation) and of Bona, Sun and Zhang \cite{bsz2002}. The first approach, which was later improved further for KdV and also adapted for NLS on the half-line by Holmer \cite{h2006, h2005}, is based on expressing the relevant forced linear IBVP as a chain of IVPs, thus allowing one to take advantage of the powerful Fourier analysis machinery but, at the same time, signifying a departure from the IBVP framework. The second approach, which has also been employed for NLS on the half-line \cite{bsz2018}, relies on solving the forced linear IBVP via a Laplace transform in the temporal variable, in contrast with the spatial (Fourier) transform used in the case of the IVP.

A novel approach was recently introduced for proving well-posedness of IBVPs for nonlinear evolution equations. This approach, which has already been  implemented  for the NLS, KdV and ``good'' Boussinesq equations on the half-line \cite{fhm2017, fhm2016, hm2015}, overcomes the lack of a proper spatial  transform in the IBVP setting by exploiting the linear solution formulae produced via the unified transform, also known as the Fokas method \cite{f1997, f2008}. 
Fokas' unified transform can be employed for solving linear evolution equations of arbitrary spatial order, supplemented with any type of admissible  boundary data (including  non-separable ones) and formulated in any number of spatial dimensions. In this light, taking also into account that no classical spatial transform exists for IBVPs involving linear evolution equations of order higher than two, the unified transform can be regarded as the analogue of the Fourier transform in the case of linear IBVPs. As such, it comes forth as the natural way of defining the iteration map to be used for showing well-posedness of nonlinear IBVPs via contraction mapping.

The essence of the new approach lies in the analysis of the pure linear IBVP, which corresponds to the case of zero initial datum and zero forcing. Indeed, the correct space for the boundary datum of the nonlinear problem is discovered in the process of estimating the unified transform solution formula of the pure linear problem. This step makes crucial use of the boundedness of the Laplace transform in $L^2$. Furthermore, it reveals the regularity of the solution of the linear IVP with respect to the boundary variable(s). In one dimension, this leads to the IVP time estimates previously obtained in \cite{kpv1991} (see also Theorem 4 in \cite{fhm2017} and Theorem~\ref{2d-nls-fls-ivp-1d-t}); in higher dimensions, it provides the extension of these estimates, as given in Theorems \ref{2d-nls-ls-ivp-t} and \ref{2d-nls-fls-ivp-t}.

Before stating our results precisely, we define the function spaces needed. The half-plane Sobolev space $H^s(\mathbb R\times \mathbb R^+)$ for the initial data is defined as a restriction of the Sobolev space $H^s(\mathbb R^2)$ by
$$
H^s(\mathbb R\times \mathbb R^+)
=
\Big\{ 
f \in L^2(\mathbb R\times \mathbb R^+): 
\no{f}_{H^s(\mathbb R\times \mathbb R^+)}
\doteq
\inf_{F|_{\mathbb R\times \mathbb R^+} = f} \no{F}_{H^s(\mathbb R^2)}<\infty
\Big\}.
$$
Furthermore, the space $B_T^s$ for the boundary data is defined as
\eee{
B_T^s 
=
\Big\{
g\in L^2(\mathbb R_{x_1}\times [0, T]): 
\no{g}_{B_T^s} 
\doteq
\no{g}_{H_{x_1}^0 H_t^{\frac{2s+1}{4}}} + \no{g}_{H_{x_1}^s H_t^{\frac{1}{4}}} <\infty
\Big\}, 
\nn
}
where the two components of the norm $\no{\cdot}_{B_T^s}$ are given by
\sss{\label{2d-nls-bst-def}
\ddd{
&\no{g}_{H_{x_1}^0 H_t^{\frac{2s+1}{4}}}
=
\left(
\int_{k_1\in\mathbb R} 
\no{e^{ik_1^2t} \, \what{g}^{x_1}(k_1, t)}_{H_t^{\frac{2s+1}{4}}(0, T)}^2 dk_1
\right)^{\frac 12},
\label{2d-nls-d1t-def}
\\
&\no{g}_{H_{x_1}^s H_t^{\frac{1}{4}}}
=
\left(
\int_{k_1\in\mathbb R} 
\left(1+k_1^2\right)^s 
\no{e^{ik_1^2t} \, \what{g}^{x_1}(k_1, t)}_{H_t^{\frac 14}(0, T)}^2 dk_1
\right)^{\frac 12},
\label{2d-nls-d2t-def}
}
}
with $\what g^{x_1}(k_1, t)$ denoting the Fourier transform of $g(x_1, t)$ with respect to $x_1$, i.e.
\eee{\label{2d-nls-ft-x1-def}
\what g^{x_1}(k_1, t)
=
\int_{x_1\in\mathbb R} e^{-ik_1x_1} g(x_1, t) dx_1.
}
In fact, as shown in Section \ref{2d-nls-ibvp-rv-s}, a straightforward manipulation of the global-in-time counterparts of the norms \eqref{2d-nls-bst-def} reveals that for $s\geqslant -\frac 12$ the space $B_T^s$ can be regarded as a restriction on $\mathbb R_{x_1}\times [0, T]$ of the space
$$
B^s 
=
X^{0, \frac{2s+1}{4}} \cap X^{s, \frac 14},
\quad
\no{g}_{B^s} 
\doteq
\no{g}_{X^{0, \frac{2s+1}{4}}}
+
\no{g}_{X^{s, \frac 14}},
$$
where $X^{s, b}$ is the usual Bourgain space associated with the NLS IVP on the line, i.e.
$$
X^{s, b}
\!=\!
\bigg\{g \in L^2(\mathbb R_{x_1}\times \mathbb R_t)\!: \no{g}_{X^{s, b}}
\doteq
\!\left(\int_{k_1\in\mathbb R}\int_{\tau\in\mathbb R}
\!\! \left(1+k_1^2\right)^s 
\!
\left(1+\left|\tau+k_1^2\right|\right)^{2b}
\left|\widehat g(k_1, \tau)\right|^2 d\tau dk_1
\right)^{\frac 12}\!\! < \infty
\bigg\}
$$
with $\what g(k_1, \tau)$ here denoting the spatiotemporal Fourier transform of $g(x_1, t)$, that is
\eee{\label{2d-nls-g-2d-ft}
\widehat g(k_1, \tau)
=
\int_{x_1\in\mathbb R}\int_{t\in\mathbb R}
e^{-ik_1 x_1-i\tau t} g(x_1, t) dtdx_1.
}

We shall show that the NLS IBVP \eqref{2d-nls-nls-ibvp} is locally well-posed  in the sense of Hadamard, namely that it possesses a unique solution which depends continuously on the prescribed initial and boundary data. 
We note that for $s>1$ the data must satisfy the compatibility condition
\eee{\label{2d-nls-comp-cond}
u_0(x_1, 0) = g_0(x_1, 0) \quad \forall x_1\in\mathbb R.
}
The precise statement of our main result is as follows.
\begin{theorem}[\b{Well-posedness of NLS on the half-plane}]
\label{2d-nls-nls-ibvp-wp-t}
Suppose $1<s \leqslant \tfrac 32$.  For initial data 
$u_0\in H^s(\mathbb R_{x_1}\times \mathbb R_{x_2}^+)$
and boundary data $g_0\in B_T^s$ satisfying the compatibility condition \eqref{2d-nls-comp-cond}, IBVP \eqref{2d-nls-nls-ibvp}   
for NLS on the half-plane has a unique solution 
$
u\in C([0, T^*]; H^s(\mathbb R_{x_1}\times \mathbb R_{x_2}^+)),
$
which admits the estimate
\eee{\nn
\sup_{t\in [0, T^*]} \no{u(t)}_{H^s(\mathbb R_{x_1}\times \mathbb R_{x_2}^+)}
\leqslant 
2c_s \no{\left(u_0, g_0\right)}_{D}, \quad c_s=c\left(s\right)>0,
}
where 
$
\no{\left(u_0, g_0\right)}_{D}
=
\no{u_0}_{H^s(\mathbb R_{x_1}\times \mathbb R_{x_2}^+)} + \no{g_0}_{B_T^s}
$
and the lifespan $T^*$ is given by
\eee{\nn
T^*= \min\left\{ T, \, c_{s, p}  \no{(u_0, g_0)}_{D}^{-2\left(p-1\right)}\right\}, \quad  c_{s, p}=c\left(s, p\right)>0.
}
Moreover, the data-to-solution map $\left\{u_0, g_0\right\}\mapsto u$ is locally Lipschitz continuous.
\end{theorem}

We note that the vast majority of results  in the literature on the Hadamard well-posedness of IBVPs for NLS in higher than one spatial dimensions refer to the case of zero (homogeneous) boundary data; see, for example, Brezis and Gallouet \cite{bg1980}, Y. Tsutsumi \cite{t1983}, M. Tsutsumi \cite{t1989, t1991}, and Burq, G\'erard and Tzvetkov \cite{bgt2003, bgt2004}.
Indeed, to the best of our knowledge, the only well-posedness results available for the non-homogeneous NLS IBVP \eqref{2d-nls-nls-ibvp} are those in the recent preprints by Audiard \cite{a2017} and Ran, Sun and Zhang \cite{rsz2017}. We emphasize, however, that the results of the present work are established via an entirely different  method, namely by advancing into two spatial dimensions the novel approach which was recently introduced for one-dimensional IBVPs \cite{fhm2017, fhm2016, hm2015} and which relies on Fokas' unified transform solution formulae. Taking into account the wide range of applicability of the unified transform, the approach developed in our work could be further adapted and generalized for showing well-posedness of IBVPs involving other well-known evolution equations in two as well as in higher dimensions.

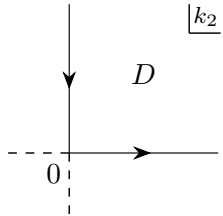
\begin{figure}[ht]
\centering
\vspace{2.2cm}
\begin{tikzpicture}[scale=1.1]
\pgflowlevelsynccm
\draw[line width=.5pt, black, dashed](0,0)--(-0.8,0);
\draw[line width=.5pt, black, dashed](0,0)--(0,-0.8);
\draw[line width=.5pt, black](1.45,1.8)--(1.45,1.45);
\draw[line width=.5pt, black](1.45,1.45)--(1.83,1.45);
\node[] at (1.59, 1.65) {\fontsize{8}{8} $k_2$};
\node[] at (-0.25, -0.25) {\fontsize{10}{10} $0$};
\draw[middlearrow={Stealth[scale=1.3, reversed]}] (0,0) -- (90:1.8);
\draw[middlearrow={Stealth[scale=1.3]}] (0,0) -- (0:1.8);
\node[] at (0.9, 0.95) {\fontsize{10}{10}\it $D$};
\end{tikzpicture}
\vspace{7mm}
\caption{The region $D$ and its positively oriented boundary $\p D$.}
\label{2d-nls-dplus}
\end{figure}

Theorem \ref{2d-nls-nls-ibvp-wp-t} will be established via a contraction mapping argument. Hence, a key role in the analysis will be played by the forced linear analogue of the nonlinear problem \eqref{2d-nls-nls-ibvp}, which reads
\sss{\label{2d-nls-fls-ibvp}
\ddd{
&i u_t + u_{x_1x_1}+u_{x_2x_2}= f \in C([0, T];  H^s(\mathbb R_{x_1}\times \mathbb R_{x_2}^+)),
\quad &&
(x_1, x_2, t)\in \mathbb R\times \mathbb R^+ \times (0, T),
\label{2d-nls-fls-eq}
\\
&u(x_1, x_2, 0) = u_0(x_1, x_2)\in H^s(\mathbb R_{x_1}\times \mathbb R_{x_2}^+), && (x_1, x_2)\in \mathbb R\times \overline{\mathbb R^+},
\\
&u(x_1, 0, t) = g_0(x_1, t) \in B_T^s, && (x_1, t)\in \mathbb R \times [0, T].
}
}
The  forced linear IBVP \eqref{2d-nls-fls-ibvp} will be estimated by taking advantage of the following explicit solution formula obtained via the unified transform of Fokas:
\ddd{\label{2d-nls-fls-ibvp-utm-sol-T}
u(x_1, x_2, t)
&= S\big[ u_0, g_0; f\big](x_1, x_2, t)
\\
&=
\frac{1}{(2\pi)^2}
\int_{k_1\in\mathbb R}
\int_{k_2\in\mathbb R}
e^{ik_1x_1+ik_2x_2-i(k_1^2+k_2^2)t}
\widehat u_0(k_1, k_2)
dk_2
dk_1
\nn\\
&
\quad
-\frac{1}{(2\pi)^2}
\int_{k_1\in\mathbb R}
\int_{k_2\in\p D}
e^{ik_1x_1+ik_2x_2-i(k_1^2+k_2^2)t}
\widehat u_0(k_1,-k_2)
dk_2
dk_1
\nn\\
&\quad
-\frac{i}{(2\pi)^2}
\int_{k_1\in\mathbb R}
\int_{k_2\in\mathbb R}
e^{ik_1x_1+ik_2x_2-i(k_1^2+k_2^2)t}
\int_{t'=0}^t 
e^{i(k_1^2+k_2^2)t'}\widehat f(k_1, k_2, t')dt'dk_2 dk_1
\nn\\
&
\quad
+\frac{i}{(2\pi)^2}
\int_{k_1\in\mathbb R}
\int_{k_2\in\p D}
e^{ik_1x_1+ik_2x_2-i(k_1^2+k_2^2)t}
\int_{t'=0}^t 
e^{i(k_1^2+k_2^2)t'}\widehat f(k_1, -k_2, t')dt'
dk_2 dk_1
\nn\\
&\quad
+\frac{1}{(2\pi)^2}
\int_{k_1\in\mathbb R}\int_{k_2\in\p D}
e^{ik_1x_1+ik_2x_2-i(k_1^2+k_2^2)t}
\, 2k_2 \, \what{g_0}^T(k_1, -k_1^2-k_2^2)
dk_2dk_1.
\nn
}
In the above formula, $\what u_0$ stands for the half-plane Fourier transform of the initial datum $u_0$, i.e.
\eee{\label{2d-nls-u0h-def}
\what u_0(k_1, k_2) 
= 
\int_{x_1\in\mathbb R} \int_{x_2\in\mathbb R^+} e^{-ik_1x_1-ik_2x_2} u_0(x_1, x_2) dx_2 dx_1,
}
which is well-defined for $k_1\in\mathbb R$ and $k_2\in \overline{\mathbb C^-} = \left\{k_2\in\mathbb C: \text{Im}(k_2)\leqslant 0\right\}$ due to the fact that $x_2\geqslant 0$.
Moreover, $\what{g_0}^T$ denotes the spatiotemporal transform
\eee{\label{2d-nls-ghatt-def}
\what{g_0}^T(k_1, \tau)
=
\int_{x_1\in\mathbb R}\int_{t=0}^T
e^{-ik_1x_1-i\tau t} g_0(x_1, t) dt dx_1,
}
which makes sense for $k_1\in\mathbb R$ and $\tau\in\mathbb C$ since integration in $t$ takes place over a bounded interval.
Finally, the contour of integration $\p D$ is the positively oriented boundary of the first quadrant $D$ of the complex $k_2$-plane, as shown in Figure \ref{2d-nls-dplus}. A concise derivation  of formula \eqref{2d-nls-fls-ibvp-utm-sol-T}, whose homogeneous version was first derived in \cite{f2002}, is provided in the appendix. 
\vspace{3mm}

Starting from the unified transform formula \eqref{2d-nls-fls-ibvp-utm-sol-T}, we shall establish the following result, which, apart from being instrumental in the proof of Theorem \ref{2d-nls-nls-ibvp-wp-t} for the well-posedness of the nonlinear problem \eqref{2d-nls-nls-ibvp}, is interesting in its own right in regard to the forced linear problem \eqref{2d-nls-fls-ibvp}.
\begin{theorem}[\b{Forced linear Schr\"odinger on the half-plane}] 
\label{2d-nls-fls-ibvp-t}
Suppose that $1< s \leqslant \frac 32$. Then, the unified transform formula \eqref{2d-nls-fls-ibvp-utm-sol-T} defines a solution $u=S\big[u_0, g_0; f\big]$ to the forced linear Schr\"odinger IBVP  \eqref{2d-nls-fls-ibvp} supplemented with the compatibility condition \eqref{2d-nls-comp-cond}, which satisfies the  estimate
\eee{\label{2d-nls-fls-ibvp-se}
\sup_{t\in[0,T]} \no{u(t)}_{H^s(\mathbb R_{x_1}\times\mathbb R_{x_2}^+)} 
\leqslant 
c_s
\Big(
\no{u_0}_{H^s(\mathbb R_{x_1}\times\mathbb R_{x_2}^+)}
+ 
\no{g_0}_{B_T^s}
+
\sqrt T \sup_{t\in[0,T]} \no{f(t)}_{H^s(\mathbb R_{x_1}\times\mathbb R_{x_2}^+)}
\Big),
}
where $c_s>0$ is a constant that depends only on $s$.
\end{theorem}

The literature on the IVP for the NLS equation in one and higher dimensions is extensive. Further results on the well-posedness, stability and blow-up behavior of NLS via hard analysis as well as integrability techniques can be found in  \cite{w1982, lpss1988, gs1993, kpv1993, m1993, cks1995, dz1995, ss1999, bgt2004b, mr2004,  v2007, ktv2009, lp2009, fh2011,  bm2017} and the references therein.
Moreover, for a thorough introduction to Fokas' unified transform method 
we refer the reader to the monograph \cite{f2008} and the review article \cite{fs2012}. In particular, we note that the IBVP on the half-line and on the interval for the cubic NLS equation has been studied via the integrable nonlinear extension of the unified transform by Fokas and collaborators in \cite{fis2005} and \cite{fi2004} respectively. Other works via similar techniques include \cite{fp2005, f2009, lf2012a, lf2012b, mf2015}.

\vskip 3mm
\noindent
\textbf{Structure of the paper.}
In Section \ref{2d-nls-ibvp-rv-s}, employing the Fokas unified transform formula \eqref{2d-nls-fls-ibvp-utm-sol-T} and exploiting the boundedness of the Laplace transform in $L^2$, we estimate the solution of the pure linear Schr\"odinger IBVP, i.e. of problem \eqref{2d-nls-fls-ibvp} with zero initial datum and zero forcing.  The resulting estimate is key to our analysis, as it reveals  the boundary data space $B_T^s$. 
In Section \ref{2d-nls-ivp-s}, we show that the $x_1t$-regularity of the solution of both the homogeneous and the forced linear Schr\"odinger on the whole line is described by $B_T^s$.
Then, in Section \ref{2d-nls-full-ibvp-s}, we combine the estimates of the previous sections  to deduce Theorem \ref{2d-nls-fls-ibvp-t} for the forced linear IBVP \eqref{2d-nls-fls-ibvp}.
This provides the central estimate needed for the proof of  Theorem \ref{2d-nls-nls-ibvp-wp-t} on the well-posedness of NLS on the half-plane, which is carried out in Section \ref{2d-nls-lwp-s} via a contraction mapping argument.
Finally, a concise derivation of the unified transform solution formula \eqref{2d-nls-fls-ibvp-utm-sol-T} is provided in the appendix.

%
%
%
%
%
%
%
%
\section{The Pure Linear IBVP}
\label{2d-nls-ibvp-rv-s}
The essence of the analysis of the forced linear IBVP \eqref{2d-nls-fls-ibvp} is captured in the simplest \textit{genuine} such problem, namely an IBVP with zero initial datum and forcing but with a \textit{nonzero boundary datum} (note that the case of zero boundary datum can be reduced to an IVP). 
In fact, this ``model problem'' can be simplified further by assuming a boundary datum with \textit{compact support} in the $t$-variable. Thus, we are led to the following problem, which we identify as the \textit{pure linear IBVP}:
\sss{\label{2d-nls-ibvp-rv}
\ddd{
&iv_t+v_{x_1x_1}+v_{x_2x_2} = 0, \quad  &&(x_1, x_2, t)\in \mathbb R\times \mathbb R^+ \times (0, 2),
\label{2d-nls-ibvp-rv-eq} 
\\
&v(x_1, x_2, 0)= 0, && (x_1, x_2)\in \mathbb R\times \overline{\mathbb R^+},  
\label{2d-nls-ibvp-rv-ic} 
\\
&v(x_1, 0, t) = g(x_1, t), &&(x_1, t)\in \mathbb R \times [0, 2], 
\label{2d-nls-ibvp-rv-bc}
\\
&\text{supp}(g)\subset \mathbb R_{x_1}\times (0, 2).
\label{2d-nls-cs}
}
}

As we shall see below, the estimation of the solution of the pure linear IBVP \eqref{2d-nls-ibvp-rv} in the usual Sobolev space $H^s(\mathbb R\times \mathbb R^+)$ reveals the ``exotic'' Bourgain-type space $B^s$ and, in turn, the restriction space $B_T^s$ for the boundary datum of the NLS IBVP \eqref{2d-nls-nls-ibvp}.
We note that, contrary to the half-line, where placing the initial datum in $H^s(\mathbb R^+)$ automatically requires a Sobolev boundary datum (see Theorem 5 in \cite{fhm2017}), in the case of the half-plane it is not \textit{a priori} clear what is the appropriate space for the boundary datum.
In this respect, the pure linear IBVP holds a central role in the analysis of both the linear and the nonlinear problem. 
Another source of motivation for the boundary data space $B_T^s$ is the linear Schr\"odinger IVP (see Theorem \ref{2d-nls-ls-ivp-t}),  in a less direct way, however.

The general Fokas unified transform formula \eqref{2d-nls-fls-ibvp-utm-sol-T} yields the solution of IBVP \eqref{2d-nls-ibvp-rv} as
\ddd{
v(x_1, x_2, t)
&= S\big[0, g; 0\big](x_1, x_2, t)
\nn\\
&=
\frac{1}{(2\pi)^2}
\int_{k_1\in\mathbb R}\int_{k_2\in\p D}
e^{ik_1x_1+ik_2x_2-i(k_1^2+k_2^2)t}
\, 2k_2 \, \what g(k_1, -k_1^2-k_2^2)
dk_2dk_1,
\label{2d-nls-fls-ibvp-rv-utm-sol-T}
}
where the contour $\p D$ is depicted in Figure \ref{2d-nls-dplus}. Importantly, note that thanks to the support assumption \eqref{2d-nls-cs} the boundary datum $g$ appears in the above formula through its two-dimensional \textit{Fourier transform} \eqref{2d-nls-g-2d-ft} instead of the ``truncated'' transform \eqref{2d-nls-ghatt-def}. 
\begin{theorem}[\b{Pure IBVP estimate}]
\label{2d-nls-ibvp-rv-t}
The solution $v=S\big[0, g; 0\big]$ of the pure linear IBVP \eqref{2d-nls-ibvp-rv} given by formula \eqref{2d-nls-fls-ibvp-rv-utm-sol-T} admits the estimate
\eee{
\sup_{t\in[0, 2]} \no{v(t)}_{H^s(\mathbb R_{x_1}\times \mathbb R_{x_2}^+)}
\leqslant 
c_s \no{g}_{B^s}, \quad s\geqslant 0, \label{2d-nls-ibvp-rv-se}
}
\end{theorem}

\begin{remark}\label{2d-nls-space-r}
In fact, one can also prove the following estimate:
\eee{\label{2d-nls-ibvp-rv-te}
\sup_{x_2\in \overline{\mathbb R^+}} \no{v(x_2)}_{B^s}
\leqslant 
c_s \no{g}_{B^s},  \quad s\in \mathbb R,
}
which allows one to carry out the contraction for the nonlinear IBVP in the space $C([0, T^*]; H^s(\mathbb  R))\cap C(\mathbb R_{x_2}^+; B_{T^*}^s)$ instead of the Hadamard space of Theorem \ref{2d-nls-nls-ibvp-wp-t}.
\end{remark}

\begin{Proof}[Proof of Theorem \ref{2d-nls-ibvp-rv-t}]
We decompose formula \eqref{2d-nls-fls-ibvp-rv-utm-sol-T} in two parts as $v = v_1 + v_2$ where  
\eee{
v_1(x_1, x_2, t)
=
\frac{1}{(2\pi)^2}
\int_{k_1\in\mathbb R}
\int_{k_2=0}^\infty
e^{ik_1x_1-k_2x_2-i(k_1^2-k_2^2)t}
\, 2k_2\widehat g(k_1, -k_1^2+k_2^2) 
dk_2 dk_1
\label{2d-v1-def}
}
corresponds to the  imaginary axis portion of $\p D$ and 
\eee{
v_2(x_1, x_2, t)
=
\frac{1}{(2\pi)^2}
\int_{k_1\in\mathbb R}\int_{k_2=0}^\infty
e^{ik_1 x_1+ik_2x_2-i(k_1^2+k_2^2)t}
\, 2k_2\widehat g (k_1, -k_1^2-k_2^2)
dk_2 dk_1
\label{2d-v2-def}
}
corresponds to the real axis part of $\p D$. The estimation of $v_2$ is easier, while that of $v_1$ is more challenging and relies crucially on the boundedness of the Laplace transform in $L^2$.
\vskip 3mm
\noindent
\textit{\b{Estimation along the real axis.}}
Let $V_2(x_1, x_2, t)$ be a global-in-space function  defined via the two-dimensional Fourier transform  
\eee{
\widehat V_2(k_1, k_2, t)
=
\left\{
\arraycolsep=4pt
\def\arraystretch{1.2}
\begin{array}{ll}
e^{-i(k_1^2+k_2^2)t} \, 2k_2\widehat g(k_1, -k_1^2 -k_2^2), 
&k_2>0,
\\
0, & k_2\leqslant 0,
\end{array}
\right.
\nn
}
so that $V_2\big|_{x_2\in\mathbb R^+} = v_2$.
In turn, we have
$$
\no{v_2(t)}_{H^s(\mathbb R_{x_1}\times \mathbb R_{x_2}^+)}^2
\leqslant
\no{V_2(t)}_{H^s(\mathbb R^2)}^2
=
\int_{k_1\in\mathbb R}\int_{k_2\in \mathbb R}
\left(1+k_1^2+k_2^2\right)^s 
\big|\what V_2(k_1, k_2, t)\big|^2 dk_2 dk_1
$$
and hence for $s\geqslant 0$  we find
\ddd{
\no{v_2(t)}_{H^s(\mathbb R_{x_1}\times \mathbb R_{x_2}^+)}^2
&\lesssim
\int_{k_1\in\mathbb R}\int_{k_2=0}^\infty
\left(1+k_1^2\right)^s k_2^2
\left|\widehat g (k_1, -k_1^2-k_2^2)\right|^2 dk_2 dk_1
\nn\\
&\quad
+
\int_{k_1\in\mathbb R}\int_{k_2=0}^\infty
\left(k_2^2\right)^{s+1}
\left|\widehat g (k_1, -k_1^2-k_2^2)\right|^2 dk_2 dk_1.
\nn
}
Making the change of variable $k_2 = \sqrt{-\tau}$ and using the property ${\what \varphi}^t(\tau-a)=\wwhat{e^{iat}\varphi(t)}^t(\tau)$, the above inequality becomes
\ddd{
\no{v_2(t)}_{H^s(\mathbb R_{x_1}\times \mathbb R_{x_2}^+)}^2
&\lesssim
\int_{k_1\in\mathbb R}\int_{\tau=-\infty}^0
\left(1+k_1^2\right)^s |\tau|^{\frac 12}
\Big|
\wwhat{e^{ik_1^2t} \what g^{x_1}(k_1, t)}^t(\tau)
\Big|^2 d\tau dk_1
\nn\\
&\quad
+
\int_{k_1\in\mathbb R}\int_{\tau=-\infty}^0
|\tau|^{s+\frac 12}
\Big| 
\wwhat{e^{ik_1^2t} \what g^{x_1}(k_1, t)}^t(\tau)
\Big|^2 d\tau dk_1.
\nn
}
Thus, since $s\geqslant 0 > -\frac 12$, we have
\sss{\label{2d-nls-bs-12}
\ddd{
\no{v_2(t)}_{H^s(\mathbb R_{x_1}\times \mathbb R_{x_2}^+)}^2
&\lesssim
\int_{k_1\in\mathbb R}
\left(1+k_1^2\right)^s
\no{e^{ik_1^2t} \what g^{x_1}(k_1, t)}_{H^{\frac{1}{4}}(\mathbb R_t)}^2 dk_1
\label{2d-nls-bs-1}
\\
&\quad
+
\int_{k_1\in\mathbb R}
\no{e^{ik_1^2t} \what g^{x_1}(k_1, t)}_{H^{\frac{2s+1}{4}}(\mathbb R_t)}^2 dk_1.
\label{2d-nls-bs-2}
}
}
The two terms on the right-hand side of the above estimate are the global-in-time counterparts of the components \eqref{2d-nls-bst-def} of the $B_T^s$-norm. In this connection, we note that
\sss{\label{2d-nls-bs-eq-norms}
\ddd{
\no{g}_{X^{0, \frac{2s+1}{4}}}^2
&=
\int_{k_1\in\mathbb R} 
\int_{\tau\in\mathbb R}
\left(1+\left|\tau+k_1^2\right|\right)^{2\cdot \frac{2s+1}{4}} 
\left|\what g(k_1, \tau)\right|^2 d\tau dk_1
\nn\\
&\simeq
\int_{k_1\in\mathbb R} 
\int_{\tau\in\mathbb R}
\left(1+\tau^2\right)^{\frac{2s+1}{4}} 
\left|\what g(k_1, \tau-k_1^2)\right|^2 d\tau dk_1
=
\eqref{2d-nls-bs-2}
}
and, similarly, 
\eee{
\no{g}_{X^{s, \frac{1}{4}}}^2
\simeq
\int_{k_1\in\mathbb R} 
\left(1+k_1^2\right)^s
\int_{\tau\in\mathbb R}
\left(1+\tau^2\right)^{\frac{1}{4}} 
\left|\what g(k_1, \tau-k_1^2)\right|^2 d\tau dk_1
=
\eqref{2d-nls-bs-1}.
}
}
Thus, by the definition of the $B^s$-norm, estimate \eqref{2d-nls-bs-12} is equivalent to
\eee{\label{2d-nls-v2-se} 
\no{v_2(t)}_{H^s(\mathbb R_{x_1}\times \mathbb R_{x_2}^+)}
\lesssim
\no{g}_{B^s}.
}
\vskip 3mm
\noindent
\textit{\b{Estimation along the imaginary axis.}}
We employ the physical space equivalent Sobolev norm
\eee{\label{2d-nls-hs-l2-frac}
\no{v_1(t)}_{H^s(\mathbb R_{x_1}\times \mathbb R_{x_2}^+)}^2
=
\sum_{|\alpha|\leqslant \left\lfloor s\right\rfloor}
\no{\p_x^\alpha v_1(t)}_{L^2(\mathbb R_{x_1}\times \mathbb R_{x_2}^+)}^2
+
\sum_{|\alpha|=\left\lfloor s\right\rfloor}
\no{\p_x^{\alpha} v_1(t)}_\beta^2, \quad s\geqslant 0,
}
where  $x=(x_1, x_2)$ and $\p_x^\alpha = \p_{x_1}^{\alpha_1}\p_{x_2}^{\alpha_2}$ with $|\alpha|=\alpha_1+\alpha_2$, and for $\beta = s-\left\lfloor s \right\rfloor\in (0, 1)$ we define
\ddd{
\no{\p_x^{\alpha} v_1(t)}_\beta^2
&=
\iint_{x\in\mathbb R\times \mathbb R^+} 
\iint_{y\in\mathbb R\times \mathbb R^+} 
\frac{\left|\p_x^{\alpha} v_1(x, t) - \p_x^{\alpha} v_1(y, t) \right|^2}{\left|x-y\right|^{2(1+\beta)}}
\, dy dx
\nn\\
&\simeq
\int_{x\in\mathbb R\times \mathbb R^+}\int_{z\in\mathbb R\times \mathbb R^+} \frac{\left| \p_x^{\alpha} v_1(x+z, t)- \p_x^{\alpha} v_1(x, t) \right|^2}{\left|z\right|^{2(1+\beta)}}\, dz dx.
\label{2d-frac-norm-def}
}
\vskip 3mm
\noindent
\textbf{Estimation of the integer part.}
We begin with the estimation of $\no{\p_x^\alpha v_1(t)}_{L^2(\mathbb R_{x_1}\times \mathbb R_{x_2}^+)}$ for all $|\alpha|\in\mathbb N_0$ with $|\alpha|\leqslant \left\lfloor s \right\rfloor$. 
Differentiating formula \eqref{2d-v1-def}, we have
$$
\p_x^\alpha v_1(x_1, x_2, t) 
\simeq
\int_{k_2=0}^\infty e^{-k_2x_2} G(k_2, x_1, t) dk_2,
$$
where 
\eee{\label{2d-nls-Q-def}
G(k_2, x_1, t)
=
\int_{k_1\in\mathbb R} e^{ik_1x_1-i(k_1^2-k_2^2)t}
k_1^{\alpha_1} k_2^{\alpha_2} k_2\widehat g(k_1, -k_1^2+k_2^2) dk_1.
}
Thus,  
\eee{
\no{\p_x^\alpha v_1(t)}_{L^2(\mathbb R_{x_1}\times \mathbb R_{x_2}^+)}^2
\simeq
\int_{x_1\in\mathbb R}
\big\|\mathcal L\left\{G(k_2, x_1, t)\right\} (x_2)\big\|_{L^2(\mathbb R_{x_2}^+)}^2 dx_1,
\nn
}
where $\mathcal L\{G\}$ denotes the Laplace transform of $G$ with respect to $k_2$, i.e. 
$$
\mathcal L\left\{G(k_2, x_1, t)\right\} (x_2)
=
\int_{k_2=0}^\infty e^{-k_2x_2} G(k_2, x_1, t) dk_2.
$$
\begin{lemma}[\b{$L^2$-boundedness of the Laplace transform}]
\label{2d-nls-lap-l}
The map
$$
\mathcal L: \phi  \mapsto  \int_{k_2=0}^\infty e^{-k_2x_2} \phi(k_2) dk_2
$$
is bounded from $L_{k_2}^2(0,\infty)$ into $L_{x_2}^2(0, \infty)$ with
$$
\no{\mathcal L\left\{\phi\right\}}_{L_{x_2}^2(0,\infty)}
\leqslant
\sqrt{\pi} \no{\phi}_{L_{k_2}^2(0,\infty)}.
$$
\end{lemma}
A proof of Lemma \ref{2d-nls-lap-l} is available in \cite{fhm2017}. Using the relevant estimate for $\phi(k_2)=G(k_2, x_1, t)$, we obtain
\ddd{
\no{\p_x^\alpha v_1(t)}_{L^2(\mathbb R_{x_1}\times \mathbb R_{x_2}^+)}^2
&\lesssim
\int_{x_1 \in \mathbb R}
\int_{k_2=0}^\infty
\bigg|
\int_{k_1\in\mathbb R} 
e^{ik_1x_1-i(k_1^2-k_2^2)t}
k_1^{\alpha_1} k_2^{\alpha_2} k_2\widehat g(k_1, -k_1^2+k_2^2)  dk_1
\bigg|^2
dk_2
dx_1
\nn\\
&\simeq
\int_{k_2=0}^\infty
\int_{k_1\in\mathbb R} 
\left(k_1^2\right)^{\alpha_1} \left(k_2^2\right)^{\alpha_2}
k_2^2 
\left|
\widehat g(k_1, -k_1^2+k_2^2)  
\right|^2
dk_1 dk_2
\nn
}
after also applying Parseval-Plancherel in $x_1$ and $k_1$. 

Therefore, 
\ddd{
\sum_{|\alpha|\leqslant \left\lfloor s \right \rfloor} 
\no{\p_x^\alpha v_1(t)}_{L^2(\mathbb R_{x_1}\times \mathbb R_{x_2}^+)}^2
&\lesssim
\sum_{|\alpha|\leqslant \left\lfloor s \right \rfloor} 
\int_{k_2=0}^\infty
\int_{k_1\in\mathbb R} 
\left(k_1^2\right)^{\alpha_1} \left(k_2^2\right)^{\alpha_2}
k_2^2 
\left|
\widehat g(k_1, -k_1^2+k_2^2)  
\right|^2
dk_1 dk_2
\nn\\
&\simeq
\int_{k_2=0}^\infty
\int_{k_1\in\mathbb R} 
\Bigg(
\sum_{|\alpha|=0}^{\left\lfloor s \right \rfloor}
\left(k_1^2 + k_2^2\right)^{|\alpha|}
\Bigg)
k_2^2
\left|
\widehat g(k_1, -k_1^2+k_2^2)  
\right|^2
dk_1 dk_2
\nn\\
&\lesssim
\int_{k_2=0}^\infty
\int_{k_1\in\mathbb R} 
\left(k_1^2 + k_2^2\right)^s
k_2^2
\left|
\widehat g(k_1, -k_1^2+k_2^2)  
\right|^2
dk_1 dk_2.
\label{2d-nls-v1-temp}
}
Then, since 
$\left(k_1^2 + k_2^2\right)^s
\lesssim 
\left(1+k_1^2\right)^s+\left(k_2^2\right)^s
$
for $s\geqslant 0$,
we infer
\ddd{
\sum_{|\alpha|\leqslant \left\lfloor s \right \rfloor} 
\no{\p_x^\alpha v_1(t)}_{L^2(\mathbb R_{x_1}\times \mathbb R_{x_2}^+)}^2
&\lesssim
\int_{k_2=0}^\infty
\int_{k_1\in\mathbb R} 
\left(1+k_1^2\right)^s
k_2^2
\left|
\widehat g(k_1, -k_1^2+k_2^2)  
\right|^2
dk_1 dk_2
\nn\\
&\quad
+
\int_{k_2=0}^\infty
\int_{k_1\in\mathbb R} 
\left(k_2^2\right)^s k_2^2
\left|
\widehat g(k_1, -k_1^2+k_2^2)  
\right|^2
dk_1 dk_2.
\nn
}
Finally, making the change of variable $k_2=\sqrt \tau$,  we conclude in view of \eqref{2d-nls-bs-eq-norms} that 
\ddd{
\sum_{|\alpha|\leqslant \left\lfloor s \right \rfloor} 
\no{\p_x^\alpha v_1(t)}_{L^2(\mathbb R_{x_1}\times \mathbb R_{x_2}^+)}^2
&\lesssim
\int_{\tau=0}^\infty
\int_{k_1\in\mathbb R} 
\left(1+k_1^2\right)^s
\sqrt \tau
\left|
\widehat g(k_1, \tau-k_1^2)  
\right|^2
dk_1 d\tau
\nn\\
&\quad
+
\int_{\tau=0}^\infty
\int_{k_1\in\mathbb R} 
\tau^{s+\frac 12}
\left|
\widehat g(k_1, \tau-k_1^2)  
\right|^2
dk_1 d\tau
\nn\\
&\leqslant
\no{g}_{X^{s, \frac{1}{4}}}^2
+
\no{g}_{X^{0, \frac{2s+1}{4}}}^2
\leqslant
\no{g}_{B^s}^2.
\label{2d-nls-v1-int-est}
}
\vskip 3mm
\noindent
\textbf{Estimation of the fractional part.}
We  proceed to the estimation of  $\no{\p_x^\alpha v_1(t)}_\beta$ where now  
$|\alpha|= \alpha_1+\alpha_2=\left\lfloor s \right\rfloor\in\mathbb N_0$.  
By the definition \eqref{2d-frac-norm-def} of the fractional norm and Parseval-Plancherel  in $x_1$ and $k_1$, we have
\ddd{
\no{\p_x^\alpha v_1(t)}_\beta^2
&\simeq
\int_{z_1\in \mathbb R}\int_{x_2=0}^\infty \int_{z_2=0}^\infty 
\frac{1}{\left|z\right|^{2(1+\beta)}} 
\int_{x_1\in\mathbb R}
\Bigg|
\int_{k_1\in\mathbb R}
e^{ik_1x_1} 
\bigg[
k_1^{\alpha_1} e^{-ik_1^2 t}
\nn\\
&\quad
\cdot 
\int_{k_2=0}^\infty
e^{- k_2x_2+ik_2^2t}
\big(
e^{ik_1 z_1 - k_2 z_2} - 1
\big)  k_2^{\alpha_2} 
k_2 \what g(k_1, -k_1^2+k_2^2)
dk_2
\bigg] dk_1
\Bigg|^2
dx_1 dz_2 dx_2  dz_1
\nn\\
&\simeq
\int_{k_1\in\mathbb R} k_1^{2\alpha_1}
\int_{z_1\in \mathbb R} 
\int_{z_2=0}^\infty 
\frac{1}{\left|z\right|^{2(1+\beta)}} 
\nn\\
&\quad
\cdot
\int_{x_2=0}^\infty
\bigg|
\int_{k_2=0}^\infty
e^{- k_2x_2+ik_2^2t}
\big(
e^{ik_1 z_1 - k_2 z_2} - 1
\big)
k_2^{\alpha_2} k_2 \what g(k_1, -k_1^2+k_2^2) 
dk_2
\bigg|^2
dx_2
dz_2  dz_1 dk_1.
\nn
}
Thus, using the Laplace transform Lemma \ref{2d-nls-lap-l} in $x_2$ and $k_2$, we obtain
\eee{
\no{\p_x^\alpha v_1(t)}_\beta^2
\lesssim
\int_{k_1\in\mathbb R} k_1^{2\alpha_1}
\int_{k_2=0}^\infty
I(k_1, k_2, \beta)
k_2^{2\alpha_2} 
k_2^2\left|\what g(k_1, -k_1^2+k_2^2) \right|^2
dk_2 dk_1,
\label{2d-i-integral}
}
where 
\ddd{\label{2d-nls-I-def}
I(k_1, k_2, \beta)
&=
\int_{z_1\in \mathbb R}\int_{z_2=0}^\infty 
\frac{\left|
e^{ik_1 z_1 - k_2 z_2}
-
1
\right|^2
}
{\left(z_1^2+z_2^2\right)^{1+\beta}}
dz_2dz_1
\nn\\
&=
\frac{1}{|k_1| k_2}
\int_{\zeta_1\in \mathbb R}\int_{\zeta_2=0}^\infty 
\frac{\left|e^{i\zeta_1 - \zeta_2}-1\right|^2}
{\left(\frac{\zeta_1^2}{k_1^2}+\frac{\zeta_2^2}{k_2^2}\right)^{1+\beta}}\,
d\zeta_2d\zeta_1.
}
\begin{lemma}
\label{2d-nls-I-est-l}
The integral $I$ defined by \eqref{2d-nls-I-def} admits the bound
\eee{\label{2d-nls-I-est}
I(k_1, k_2, \beta)
\lesssim
\left(k_1^2+k_2^2\right)^\beta, \quad \beta\in (0, 1).
}
\end{lemma}

Lemma \ref{2d-nls-I-est-l} is proven after the end of the current proof. Combining  estimate \eqref{2d-nls-I-est} with inequality \eqref{2d-i-integral}, we deduce
\eee{ 
\no{\p_x^\alpha v_1(t)}_\beta^2
\lesssim
\int_{k_1\in\mathbb R} k_1^{2\alpha_1}
\int_{k_2=0}^\infty  \left(k_1^2+k_2^2\right)^\beta 
k_2^{2\alpha_2} 
k_2^2\left|\what g(k_1, -k_1^2+k_2^2) \right|^2 dk_2  dk_1.
\nn
}
Therefore, recalling that $\alpha_2=|\alpha|-\alpha_1$, we obtain
\ddd{
\sum_{|\alpha|=\left\lfloor s \right \rfloor}
\no{\p_x^\alpha v_1(t)}_\beta^2
&\lesssim
\sum_{|\alpha|=\left\lfloor s \right \rfloor}
\int_{k_1\in\mathbb R} 
\int_{k_2=0}^\infty  \left(k_1^2+k_2^2\right)^\beta 
k_1^{2\alpha_1} k_2^{2\alpha_2} 
k_2^2\left|\what g(k_1, -k_1^2+k_2^2) \right|^2 dk_2  dk_1
\nn\\
&\simeq
\int_{k_1\in\mathbb R} 
\int_{k_2=0}^\infty  \left(k_1^2+k_2^2\right)^\beta 
\left(k_1^2+k_2^2\right)^{\left\lfloor s \right \rfloor}
k_2^2\left|\what g(k_1, -k_1^2+k_2^2) \right|^2 dk_2  dk_1
\nn\\
&=
\int_{k_1\in\mathbb R} 
\int_{k_2=0}^\infty  \left(k_1^2+k_2^2\right)^s
k_2^2\left|\what g(k_1, -k_1^2+k_2^2) \right|^2 dk_2  dk_1.
\nn
}
The right-hand side above is the same with that of \eqref{2d-nls-v1-temp}. Hence,  proceeding as before we conclude that
\eee{\label{2d-nls-v1-frac-est}
\sum_{|\alpha|=\left\lfloor s \right \rfloor}
\no{\p_x^\alpha v_1(t)}_\beta
\lesssim
\no{g}_{B^s}, \quad s\geqslant 0.
}

The integer estimate \eqref{2d-nls-v1-int-est} and the fractional estimate \eqref{2d-nls-v1-frac-est} combined with the definition \eqref{2d-nls-hs-l2-frac} of the Sobolev norm yield 
\eee{\label{2d-nls-v1-se}
\no{v_1(t)}_{H^s(\mathbb R_{x_1}\times \mathbb R_{x_2}^+)}
\lesssim
\no{g}_{B^s}, \quad s\geqslant 0.
}
This estimate together with estimate \eqref{2d-nls-v2-se}  for $v_2$ implies estimate \eqref{2d-nls-ibvp-rv-se}  for the solution $v$ of the pure linear IBVP.

The proof of Theorem \ref{2d-nls-ibvp-rv-t} is complete.
\end{Proof}

\noindent
\textbf{Proof of Lemma \ref{2d-nls-I-est-l}.} 
Note that $I$ is even in $k_1$. Hence, without loss of generality, we estimate $I$ for $(k_1, k_2)\in \mathbb R^+\times \mathbb R^+$. There are two cases to consider: (i) $k_2=\lambda k_1$, $0\leqslant \lambda \leqslant 1$, and (ii) $k_1=\lambda k_2$, $0\leqslant \lambda \leqslant 1$.
In the first case, starting from \eqref{2d-nls-I-def} we have
\ddd{
I = I(k_1, \lambda k_1, \beta)
&=
\frac{(k_1^2)^{1+\beta}}{\lambda k_1^2}
\int_{\zeta_1\in \mathbb R}\int_{\zeta_2=0}^\infty 
\frac{\left|e^{i\zeta_1 - \zeta_2}-1\right|^2}
{\left(\zeta_1^2+\frac{\zeta_2^2}{\lambda^2}\right)^{1+\beta}}\,
d\zeta_2d\zeta_1
\nn\\
&=
\left(k_1^2\right)^{\beta} 
\int_{\zeta_1\in \mathbb R}\int_{\zeta_2=0}^\infty 
\frac{\left|
e^{i\zeta_1 - \lambda \zeta_2}
-
1
\right|^2
}
{\left(\zeta_1^2+\zeta_2^2\right)^{1+\beta}}\,
d\zeta_2d\zeta_1
=
\left(k_1^2\right)^{\beta} J_1(\lambda, \beta),
\nn
}
after  letting $\zeta_2\mapsto \lambda \zeta_2$.
Switching to polar coordinates $\zeta_1=r\cos\theta$, $\zeta_2=r\sin\theta$ yields
$$
J_1(\lambda, \beta)
=
\int_{\theta=0}^\pi \int_{r=0}^\infty 
\frac{\left|e^{\left(i\cos \theta - \lambda \sin\theta\right)r}-1\right|^2}
{r^{1+2\beta}}\,
drd\theta.
$$
For $r\ll 1$, we have
$
e^{\left(i\cos \theta - \lambda \sin\theta\right)r}
=
1 + \left(i\cos \theta - \lambda \sin\theta\right)r + O(r^2)
$
thus the integrand of $J_1$ becomes
$$
\frac{\left|e^{\left(i\cos \theta - \lambda \sin\theta\right)r}-1\right|^2}
{r^{1+2\beta}}
=
\frac{1}{r^{2\beta-1}}  \left|\left(i\cos \theta - \lambda \sin\theta\right) + O(r)\right|^2,
$$
which is integrable at $r=0$ since $\beta<1$. 
Furthermore, for $r\gg 1$ the integrand of $J_1$ becomes
$$
\frac{\left|e^{\left(i\cos \theta - \lambda \sin\theta\right)r}-1\right|^2}
{r^{1+2\beta}}
\leqslant 
\frac{(1+1)^2}{r^{1+2\beta}},
$$
which is integrable at $r=\infty$ since $\beta>0$. Thus, $J_1(\lambda, \beta) = c_{\lambda, \beta} < \infty$ and hence
\eee{\label{2d-nls-i-c1-est}
I(k_1, \lambda k_1, \beta) \lesssim \left(k_1^2\right)^\beta, \quad 0\leqslant \lambda \leqslant 1.
}
A similar argument shows that 
\eee{\label{2d-nls-i-c2-est}
I(\lambda k_2, k_2, \beta) \lesssim \left(k_2^2\right)^\beta, \quad 0\leqslant \lambda \leqslant 1.
}
Combining estimates \eqref{2d-nls-i-c1-est} and \eqref{2d-nls-i-c2-est} concludes the proof of Lemma \ref{2d-nls-I-est-l}. \hfill $\blacksquare$

%
%
%
%
%
%
%
%
\section{Linear IVP Estimates}
\label{2d-nls-ivp-s}

The pure linear IBVP \eqref{2d-nls-ibvp-rv} will be combined with appropriate linear IVPs to yield  the forced linear IBVP \eqref{2d-nls-fls-ibvp} via the superposition principle. The details of this (de)composition are given in Section \ref{2d-nls-full-ibvp-s}. 
Hence,  Theorem \ref{2d-nls-ibvp-rv-t} for the forced linear IBVP will be established by combining the estimate of Theorem \ref{2d-nls-ibvp-rv-t} for the pure linear IBVP with suitable estimates for the aforementioned IVPs. 
These estimates are derived below.

\vskip 3mm
\noindent
\textbf{Homogeneous linear IVP estimates.}
We begin with the  linear Schr\"odinger IVP
\sss{\label{2d-nls-ls-ivp}
\ddd{
&i U _t+  U _{x_1x_1} + U_{x_2x_2}=0,  && (x_1, x_2, t) \in \mathbb R \times \mathbb R \times (0, T),
\label{2d-nls-lsivp-eq}
\\
&U (x_1, x_2, 0)= U_0(x_1, x_2)\in H^s(\mathbb R_{x_1}\times \mathbb R_{x_2}), \quad && (x_1, x_2) \in \mathbb R \times \mathbb R, \label{2d-nls-lsivp-ic}
}
}
whose solution is given by
\eee{
\label{2d-nls-ls-ivp-sol}
U (x_1, x_2, t) 
=
S\big[U_0; 0\big] (x_1, x_2, t) 
= 
\frac{1}{(2\pi)^2} 
\int_{k_1\in \mathbb R}\int_{k_2\in \mathbb R} e^{ik_1 x_1+ik_2x_2-i(k_1^2+k_2^2)t}\, \widehat U_0(k_1, k_2) dk_2 dk_1,
}
where $\widehat U_0$ is the Fourier transform of $U_0$ on the plane defined by %
\eee{\label{2d-nls-ft-def}
\what U_0(k_1, k_2)
=
\int_{x_1\in\mathbb R} \int_{x_2\in\mathbb R}
e^{-ik_1x_1-ik_2x_2}  U_0(x_1, x_2) dx_2 dx_1.
}

We shall discover below that regularity of the solution of IVP \eqref{2d-nls-ls-ivp} in the variables $x_1$ and $t$ is described by the boundary data space $B_T^s$.  This result can be regarded as the two-dimensional analogue of the time estimate of \cite{kpv1991}, which states that  the solution of the linear Schr\"odinger IVP on the line with initial data in $H^s(\mathbb R_x)$ belongs to $H^{\frac{2s+1}{4}}(0, T)$ as a function of $t$. Hence, both in one and in two dimensions, the regularity of the linear IVP with respect to the boundary variables is described by the boundary data space of the associated Dirichlet IBVP. The precise statement of our result is the following.

\begin{theorem}[\b{Homogeneous IVP estimates}]
\label{2d-nls-ls-ivp-t}
The solution $U=S\big[U_0; 0\big]$ of the homogeneous linear IVP \eqref{2d-nls-ls-ivp} given by formula \eqref{2d-nls-ls-ivp-sol} satisfies the estimates
\ddd{
&\sup_{t\in [0, T]}
\no{U(t)}_{H^s(\mathbb R_{x_1}\times \mathbb R_{x_2})} 
=
\no{U_0}_{H^s(\mathbb R_{x_1} \times \mathbb R_{x_2})}, \quad && s\in\mathbb R, 
\label{2d-nls-ls-ivp-se}
\\
&\sup_{x_2\in\mathbb R}
\no{U(x_2)}_{B_T^s} 
\leqslant 
c_s \no{U_0}_{H^s(\mathbb R_{x_1} \times \mathbb R_{x_2})}, 
&& s\geqslant 0.
\label{2d-nls-ls-ivp-te}
}
\end{theorem}

\begin{remark}
Recall that $B_T^s$ emerges naturally in the proof of estimate \eqref{2d-nls-ibvp-rv-se} as the boundary data space that allows for the solution of the pure linear IBVP to belong in $H^s(\mathbb R\times \mathbb R^+)$.
Estimate \eqref{2d-nls-ls-ivp-te} indicates an alternative, reverse path for discovering $B_T^s$, namely by investigating the $x_1t$-regularity of the solution of the linear Schr\"odinger IVP when the initial datum belongs in  $H^s(\mathbb R\times \mathbb R)$.
Together, estimates \eqref{2d-nls-ibvp-rv-se} and \eqref{2d-nls-ls-ivp-te} can be visualized as the upper and lower half of a closed loop from the data space to the solution space and back, and hence confirm that $B_T^s$ is indeed the correct boundary data space for NLS on the half-plane in the case of Sobolev initial data.
\end{remark}

\begin{Proof}[Proof of Theorem \ref{2d-nls-ls-ivp-t}]
The isometry relation \eqref{2d-nls-ls-ivp-se} follows easily from formula \eqref{2d-nls-ls-ivp-sol}  and the definition of the Sobolev norm. 
Concerning estimate \eqref{2d-nls-ls-ivp-te}, we note that  the function 
\eee{
Q(k_1, x_2, t)
=
e^{ik_1^2 t} \what{U}^{x_1}(k_1, x_2, t) 
\nn
}
involved in the $B_T^s$-norm  of $U$ satisfies for all $k_1\in\mathbb R$ the \textit{one-dimensional} IVP
\sss{\label{2d-nls-ls-ivp-ft}
\ddd{
&i Q_t  + Q_{x_2x_2}=0, &&(x_2, t) \in \mathbb R  \times (0, T),
\label{2d-nls-lsivp-eq-ft}
\\
&Q(k_1, x_2, 0) =  \what{U_0}^{x_1} (k_1, x_2), \quad &&  x_2 \in \mathbb R.
\label{2d-nls-lsivp-ic-ft}
}
}
Thus, we have the following estimate from Theorem 4 of \cite{fhm2017}:
\eee{\label{2d-nls-j-te}
\no{Q(k_1, x_2)}_{H_t^{\frac{2s+1}{4}}(0, T)}
\lesssim
\big\|\what{U_0}^{x_1} (k_1)\big\|_{H^s(\mathbb R_{x_2})}, 
\quad s\geqslant -\tfrac 12, \ k_1, x_2\in\mathbb R.
}
In turn, for all $x_2\in\mathbb R$ and $s\geqslant -\tfrac 12$ we infer
\ddd{
\no{U(x_2)}_{B_T^s}^2
&\lesssim
\int_{k_1\in\mathbb R}
\no{\what{U_0}^{x_1} (k_1, x_2)}_{H^s(\mathbb R_{x_2})}^2
 dk_1
+ \int_{k_1\in\mathbb R}
 \left(1+k_1^2\right)^s
\no{\what{U_0}^{x_1} (k_1, x_2)}_{H^0(\mathbb R_{x_2})}^2
 dk_1
\nn\\
&=
\int_{k_1\in\mathbb R}
\int_{k_2\in\mathbb R}
\Big[
\left(1+k_1^2\right)^s 
+
\left(1+k_2^2\right)^s 
\Big]
\left|\what U_0 (k_1, k_2)\right|^2 dk_2
 dk_1,  
 \nn
}
which implies estimate \eqref{2d-nls-ls-ivp-te} upon restricting $s\geqslant 0$.
\end{Proof}

\vskip 3mm
\noindent
\textbf{Forced linear IVP estimates.}
We continue with the estimation of the forced linear Schr\"odinger on the plane.
First, we need a result for the  one-dimensional  forced IVP 
\sss{\label{2d-nls-fls-ivp-1d}
\ddd{
&i w_t +  w_{xx} = f(x, t),\quad  && (x, t) \in \mathbb R  \times \mathbb R,
\label{2d-nls-flsivp-eq-1d}
\\
&w(x, 0)=0, &&  x\in \mathbb R, \label{2d-nls-flsivp-ic-1d}
}
}
whose solution is    given by
\sss{\label{2d-nls-fls-ivp-sol-comb-1d}
\ddd{
\label{2d-nls-fls-ivp-sol-1d}
w(x, t) 
= 
S\big[0; f\big](x, t)
&=
-\frac{i}{2\pi}\int_{t'=0}^t   \int_{k\in \mathbb R}  e^{ik x - ik^2(t-t')}
\what f(k, t') dk  dt' 
\quad &&
\\
&=
-i \int_{t'=0}^t  S\big[f(\cdot, t'); 0\big](x, t-t') dt',
\label{2d-nls-fls-ivp-sol-d-1d}
}
}
where $\what f$ is the whole-line Fourier transform of $f$ defined by
\eee{
\what f(k, t) = \int_{x\in\mathbb R} e^{-ikx} f(x, t) dx
}
and $S\big[f(x, t'); 0\big]$ denotes the solution of the one-dimensional homogeneous  IVP \eqref{2d-nls-ls-ivp-ft} with initial datum $f(x, t')$.

\begin{theorem}[\b{Forced IVP time estimate in one dimension}]
\label{2d-nls-fls-ivp-1d-t}
For all $x\in\mathbb R$, the solution $w=S\big[0; f\big]$ of the one-dimensional forced linear IVP \eqref{2d-nls-fls-ivp-1d} given by formula \eqref{2d-nls-fls-ivp-sol-comb-1d} admits the bounds
\sss{\label{2d-nls-1d-fls-te}
\ddd{
\no{w(x)}_{H_t^{\frac{2s+1}{4}}(0, T)}^2
&\lesssim
 \int_{t=0}^{T}
\no{f(t)}_{H^s(\mathbb R_x)}^2 dt
\nn\\
&\quad 
+\int_{t=0}^T \int_{z=0}^{T-t}
\frac{1}{z^{\frac 32+s}}
\left[
\int_{t'=t}^{t+z} 
\no{f(t')}_{H^s(\mathbb R_x)}
dt' \right]^2 dz dt,
\quad  &&\tfrac 12 < s < \tfrac 32,
\label{2d-nls-1d-fls-te-a}
\\
\no{w(x)}_{H_t^{\frac{2s+1}{4}}(0, T)}^2
&\lesssim
\int_{t=0}^T \no{f(t)}_{H^s(\mathbb R_{x})}^2 dt, 
&& s=0, \tfrac 32.
\label{2d-nls-1d-te-b}
}
}
\end{theorem}

\begin{Proof}[Proof of Theorem \ref{2d-nls-fls-ivp-t}]
Let $m=\frac{2s+1}{4}$ and note that $0\leqslant m <1$ corresponds to $-\frac 12\leqslant s<\frac 32$. Hence, in this range of $s$ the physical space equivalent $H_t^m(0, T)$-norm reads
\eee{\label{2d-nls-frac-t-def-01}
\no{w(x)}_{H_t^m(0, T)}^2
=
\no{w(x)}_{L_t^2(0, T)}^2
+
\no{w(x)}_m^2, \quad 0 \leqslant m < 1,
}
where the fractional part of the norm is defined for $m \in (0, 1)$ by
\eee{
\no{w(x)}_m^2
=
\int_{t=0}^T \int_{z=0}^{T-t} \frac{\left|w(x, t+z)-w(x, t)\right|^2}{z^{1+2m}} \, dz dt.
}

Starting from formula \eqref{2d-nls-fls-ivp-sol-d-1d} and using  Minkowski's integral inequality and the homogeneous IVP time estimate \eqref{2d-nls-j-te}, we find
\eee{\label{2d-nls-q-te-l2}
\no{w(x)}_{L_t^2(0, T)}
\lesssim
\int_{t'=0}^T 
\no{S\big[f(\cdot, t'); 0\big] (x, t-t')}_{L_t^2(0, T)} dt'
\lesssim
\int_{t'=0}^T 
\no{f(t')}_{H^{-\frac 12}(\mathbb R_{x})} dt'.
}

Moreover,  employing the Duhamel representation \eqref{2d-nls-fls-ivp-sol-d-1d}, we have
\sss{\label{2d-nls-q-te-ab}
\ddd{
\hskip -3mm
\no{w(x)}_m^2
&\leqslant
\int_{t=0}^T \int_{z=0}^{T-t}
\frac{1}{z^{1+2m}}
\left(
\displaystyle \int_{t'=0}^{T} 
\Big|S\big[f(\cdot, t'); 0\big](x, t+z-t') - S\big[f(\cdot, t'); 0\big](x, t-t')\Big| dt' 
\right)^2   dz dt
\label{2d-nls-q-te-a}
\\
&\quad
+
\int_{t=0}^T \int_{z=0}^{T-t}
\frac{1}{z^{1+2m}}
\displaystyle \left|\int_{t'=t}^{t+z} S\big[f(\cdot, t'); 0\big](x, t+z-t') dt' \right|^2 dz dt.
\label{2d-nls-q-te-b}
}
}

The term \eqref{2d-nls-q-te-a} can be estimated via Minkowski's integral inequality in the $t'$- and $zt$-integrals, followed by the homogeneous IVP time estimate  \eqref{2d-nls-j-te}. Eventually, we find 
\eee{
\eqref{2d-nls-q-te-a}
\lesssim
\int_{t=0}^{T}
\no{f(t)}_{H^s(\mathbb R_{x})}^2 dt.
\label{2d-nls-q-te-a-est}
}

For the term \eqref{2d-nls-q-te-b}, we treat the cases  $m >\frac 12$ and $m<\frac 12$ separately. In the former case, using formula \eqref{2d-nls-fls-ivp-sol-1d} and applying Cauchy-Schwarz inequality in $k$, we have 
\ddd{
\eqref{2d-nls-q-te-b}
&\lesssim
\int_{t=0}^T \int_{z=0}^{T-t}
\frac{1}{z^{1+2m}}
\left(
\int_{k\in \mathbb R} 
\int_{t'=t}^{t+z} 
\left|
\what f (k, t')
\right|  dt'
dk 
\right)^2  dz dt
\nn\\
&\lesssim
\int_{t=0}^T \int_{z=0}^{T-t}
\frac{1}{z^{1+2m}}
\int_{k\in \mathbb R}  \left(1+k^2\right)^s
\left(\int_{t'=t}^{t+z} 
\left|
\what f (k, t')
\right|  dt' \right)^2
dk  dz dt.
\nn
}
where we have made crucial use of the fact that $m>\frac 12$ implies $s>\frac 12$ and hence $\int_{k\in \mathbb R} \frac{dk}{\left(1+k^2\right)^s} <\infty$.
Thus, Minkowski's integral inequality in $t'$ and $k$ yields
\eee{\label{2d-nls-q-te-b-est}
\eqref{2d-nls-q-te-b}
\lesssim
\int_{t=0}^T \int_{z=0}^{T-t}
\frac{1}{z^{1+2m}}
\left[
\int_{t'=t}^{t+z} 
\no{f(t')}_{H^s(\mathbb R_{x})}
dt'
\right]^2 dz dt.
}
Combining estimates \eqref{2d-nls-q-te-a-est}  and \eqref{2d-nls-q-te-b-est}, we deduce the estimate
\eee{
\no{w(x)}_m^2
\lesssim
\int_{t'=0}^{T}
\no{f(t')}_{H^s(\mathbb R_{x})}^2 dt'
+
\int_{t=0}^T \int_{z=0}^{T-t}
\frac{1}{z^{1+2m}}
\left[
\int_{t'=t}^{t+z} 
\no{f(t')}_{H^s(\mathbb R_{x})}
dt'
\right]^2 dz dt, \quad \tfrac 12<m<1,
\nn
}
which, together with the $L^2$-estimate \eqref{2d-nls-q-te-l2}, implies the bound \eqref{2d-nls-1d-fls-te-a}.

For $m<\frac 12$, applying Cauchy-Schwarz inequality in $t'$ and then interchanging the $t$- and $z$-integrals, we have
\eee{
\eqref{2d-nls-q-te-b} 
\leqslant
\int_{z=0}^T
\frac{1}{z^{2m}}
 \int_{t=0}^{T-z} 
 \int_{t'=t}^{t+z} \left|S\big[f(\cdot, t'); 0\big](x, t+z-t')\right|^2 dt' dt dz.
 \nn
}
Hence, letting $t \rightarrow t - z$ and augmenting the range of integration with respect to $t$ and $t'$, we find
\ddd{
\eqref{2d-nls-q-te-b} 
&\leqslant
\int_{z=0}^T
\frac{1}{z^{2m}}
 \int_{t=z}^T 
\int_{t'=t-z}^{t} \left|S\big[f(\cdot, t'); 0\big](x, t-t')\right|^2 dt' dt dz
\nn\\
&\leqslant
\left(\int_{z=0}^T \frac{1}{z^{2m}} \, dz\right)
\int_{t'=0}^T
\no{S\big[f(\cdot, t'); 0\big](x, t-t')}_{L_t^2(0, T)}^2 dt',
\nn
}
so, integrating in $z$ (recall that $2m<1$) and using the homogeneous IVP time estimate \eqref{2d-nls-j-te}, we obtain
\eee{\label{2d-nls-kat1}
\eqref{2d-nls-q-te-b} 
\leqslant
\frac{T^{1-2m}}{1-2m}
\int_{t=0}^T
\no{f(t)}_{H^{-\frac 12}(\mathbb R_{x})}^2 dt.
}
Estimates \eqref{2d-nls-q-te-a-est} and \eqref{2d-nls-kat1} together imply
\eee{
\no{w(x)}_m^2
\lesssim
 \int_{t=0}^{T}
\no{f(t)}_{H^s(\mathbb R_{x})}^2 dt
+
\frac{T^{1-2m}}{1-2m}
\int_{t=0}^T
\no{f(t)}_{H^{-\frac 12}(\mathbb R_{x})}^2 dt, \quad 0<m<\tfrac 12,
\nn
}
thus, recalling also the $L^2$-estimate \eqref{2d-nls-q-te-l2} and the fact that $T<1$, we deduce
\eee{
\no{w(x)}_{H_t^m(0, T)}^2
\lesssim
\frac{1}{1-2m}
\int_{t=0}^T
\no{f(t)}_{H^s(\mathbb R_{x})}^2 dt,
\quad 
0 \leqslant m < \tfrac 12,
\nn
}
which is the bound \eqref{2d-nls-1d-te-b} for $s=0$.

Finally, in order to establish \eqref{2d-nls-1d-te-b} for $s=\frac 32$ (which corresponds to $m=1$) we need to estimate
$$
\no{w(x)}_{H_t^1(0, T)}^2
=
\no{w(x)}_{L_t^2(0, T)}^2
+
\no{\p_t w(x)}_{L_t^2(0, T)}^2.
$$
The first $L^2$-norm above was estimated earlier --- see \eqref{2d-nls-q-te-l2}. %
Furthermore, differentiating the Duhamel representation \eqref{2d-nls-fls-ivp-sol-d-1d} with respect to $t$, we have
\eee{\label{2d-nls-qt-te-temp1} 
\no{\p_t w(x)}_{L_t^2(0, T)}^2
\lesssim
\no{f(x)}_{L_t^2(0, T)}^2
+
\no{\int_{t'=0}^t S\big[\p_{x}^2 f(\cdot, t'); 0\big] (x, t-t')dt'}_{L_t^2(0, T)}^2.
\nn
}
Since $s>\frac 12$ , the Sobolev embedding theorem in $x$ implies
\eee{\nn
\no{f(x)}_{L_t^2(0, T)}^2
\leqslant
\int_{t=0}^T  \no{f(t)}_{L^\infty(\mathbb R_{x})}^2 dt
\leqslant
\int_{t=0}^T  \no{f(t)}_{H^s(\mathbb R_{x})}^2 dt.
}
Moreover,  similarly to the derivation of \eqref{2d-nls-q-te-l2}, we have
\eee{\nn
\no{\int_{t'=0}^t S\big[\p_{x}^2 f(t'); 0\big] (x, t-t')dt'}_{L_t^2(0, T)}
\lesssim
\int_{t'=0}^T \no{\p_{x}^2 f(t')}_{H^{-\frac 12}(\mathbb R_{x})} dt' 
\leqslant
\int_{t'=0}^T \no{f(t')}_{H^{\frac 32}(\mathbb R_{x})} dt'.
}
Therefore, applying Cauchy-Schwarz inequality in $t'$ we obtain 
$$
\no{\p_t w(x)}_{L_t^2(0, T)}^2
\lesssim
\int_{t'=0}^T \no{f(t')}_{H^s(\mathbb R_{x})}^2 dt',
$$
which combined with estimate \eqref{2d-nls-q-te-l2} yields the bound \eqref{2d-nls-1d-te-b} for $s=\frac 32$.

The proof of Theorem \ref{2d-nls-fls-ivp-t} is complete.
\end{Proof}

The bounds of Theorem \ref{2d-nls-fls-ivp-t}  will now be employed in the estimation of the two-dimensional forced linear Schr\"odinger IVP with zero initial datum:
\sss{\label{2d-nls-fls-ivp}
\ddd{
&i W_t+  W_{x_1x_1} + W_{x_2x_2}= F(x_1, x_2, t),\quad  && (x_1, x_2, t) \in \mathbb R \times \mathbb R \times \mathbb R,
\label{2d-nls-flsivp-eq}
\\
& W(x_1, x_2, 0)=0, &&  (x_1, x_2) \in \mathbb R \times \mathbb R, \label{2d-nls-flsivp-ic}
}
}
with solution
\sss{\label{2d-nls-fls-ivp-sol-comb}
\ddd{
\label{2d-nls-fls-ivp-sol}
W(x_1, x_2, t) 
&= 
S\big[0; F\big](x_1, x_2, t)
\nn\\
&=
-\frac{i}{(2\pi)^2}\int_{t'=0}^t   \int_{k_1\in \mathbb R} \int_{k_2\in \mathbb R} e^{ik_1 x_1+ik_2x_2-i(k_1^2+k_2^2)(t-t')}
\what F^x (k_1, k_2, t') dk_2 dk_1 dt' 
\quad &&
\\
&=
-i \int_{t'=0}^t  S\big[F(\cdot, \cdot, t'); 0\big](x_1, x_2, t-t') dt',
\label{2d-nls-fls-ivp-sol-d}
}
}
where $\widehat{F}$ is the Fourier transform of $F$ on the plane defined similarly to \eqref{2d-nls-ft-def} and $S\big[F(x_1, x_2, t'); 0\big]$ denotes the solution of the homogeneous IVP \eqref{2d-nls-ls-ivp} with initial datum $F(x_1, x_2, t')$.
\begin{theorem}[\b{Forced IVP estimates}]
\label{2d-nls-fls-ivp-t}
The solution $W=S\big[0; F\big]$ of the forced linear Schr\"odinger IVP \eqref{2d-nls-fls-ivp} given by \eqref{2d-nls-fls-ivp-sol-comb} satisfies the estimates
\ddd{
&\sup_{t\in [0, T]} \no{W(t)}_{H^s(\mathbb R_{x_1}\times \mathbb R_{x_2})} 
\leqslant
T \sup_{t\in [0, T]} \no{F(t)}_{H^s(\mathbb R_{x_1} \times \mathbb R_{x_2})},  \quad &&s\in\mathbb R,
\label{2d-nls-fls-ivp-se}
\\
&\sup_{x_2\in\mathbb R}
\no{W(x_2)}_{B_T^s} 
\leqslant 
c_s 
\sqrt T \sup_{t\in [0, T]} \no{F(t)}_{ H^s(\mathbb R_{x_1} \times \mathbb R_{x_2})},
&&\tfrac 12 < s \leqslant \tfrac 32.
\label{2d-nls-fls-ivp-bst}
}
\end{theorem}

\begin{Proof}[Proof of Theorem \ref{2d-nls-fls-ivp-t}]
Estimate \eqref{2d-nls-fls-ivp-se} can be derived by employing the Duhamel representation \eqref{2d-nls-fls-ivp-sol-d} in combination with the space estimate \eqref{2d-nls-ls-ivp-se} for the homogeneous  IVP.

Concerning estimate  \eqref{2d-nls-fls-ivp-bst}, we recall that 
\sss{
\ddd{
\no{W(x_2)}_{H_{x_1}^0 H_t^{\frac{2s+1}{4}}}
&=
\left(\int_{k_1\in\mathbb R} \no{e^{ik_1^2t} \what{W}^{x_1}(k_1, x_2, t)}_{H_t^{\frac{2s+1}{4}}(0, T)}^2 dk_1\right)^{\frac 12},
\label{2d-nls-kat-0}
\\
\no{W(x_2)}_{H_{x_1}^s H_t^{\frac{1}{4}}}
&=
\left(\int_{k_1\in\mathbb R} \left(1+k_1^2\right)^s \no{e^{ik_1^2t} \what{W}^{x_1}(k_1, x_2, t)}_{H_t^{\frac{1}{4}}(0, T)}^2 dk_1\right)^{\frac 12},
\label{2d-nls-kat-1}
}
}
and note that the function
\eee{
Q(k_1, x_2, t)  = e^{ik_1^2t} \what{W}^{x_1}(k_1, x_2, t) 
\nn
}
satisfies for all $k_1\in\mathbb R$ the one-dimensional IVP 
\ddd{
&i Q_t  + Q_{x_2x_2}= R(k_1, x_2, t), \quad && (x_2, t) \in \mathbb R  \times \mathbb R,
\nn\\
&Q(k_1, x_2, 0) =  0, &&  x_2 \in \mathbb R,
\nn
}
with forcing $R$ given by
\eee{\nn
R(k_1, x_2, t) = e^{ik_1^2 t} \what F^{x_1}(k_1, x_2, t).
}
Hence, for all $k_1\in\mathbb R$ the function $Q$ admits the bounds of Theorem \ref{2d-nls-fls-ivp-t}, i.e.
\sss{\label{2d-nls-1d-fls-te-q}
\ddd{
\no{Q(k_1, x_2)}_{H_t^{\frac{2s+1}{4}}(0, T)}^2
&\lesssim
\int_{t=0}^T \int_{z=0}^{T-t}
\frac{1}{z^{\frac 32+s}}
\left[
\int_{t'=t}^{t+z} 
\no{R(k_1, t')}_{H^s(\mathbb R_{x_2})}
dt' \right]^2 dz dt
\nn\\
&\quad
+\int_{t=0}^{T}
\no{R(k_1, t)}_{H^s(\mathbb R_{x_2})}^2 dt,
\quad  &&\tfrac 12 < s < \tfrac 32,
\label{2d-nls-1d-fls-te-a}
\\
\no{Q(k_1, x_2)}_{H_t^{\frac{2s+1}{4}}(0, T)}^2
&\lesssim
\int_{t=0}^T \no{R(k_1, t)}_{H^s(\mathbb R_{x_2})}^2 dt, 
&& s=0, \tfrac 32.
\label{2d-nls-1d-fls-te-q-b}
}
}

Combining the bound \eqref{2d-nls-1d-fls-te-a} with the norm \eqref{2d-nls-kat-0} and then applying Minkowski's integral inequality in $k_1$ and $t'$, we obtain
\ddd{ 
\no{W(x_2)}_{H_{x_1}^0 H_t^{\frac{2s+1}{4}}}^2
&\lesssim
\int_{t=0}^T \int_{z=0}^{T-t}
\frac{1}{z^{\frac 32 + s}}
\left[
\int_{t'=t}^{t+z} 
\left(
\int_{k_1\in\mathbb R} 
\no{R(k_1, t')}_{H^s(\mathbb R_{x_2})}^2
dk_1
\right)^{\frac 12} dt'
\right]^2
dz dt
\nn\\
&\quad
+
\int_{t'=0}^{T}
\left(
\int_{k_1\in\mathbb R} 
\no{R(k_1, t')}_{H^s(\mathbb R_{x_2})}^2 dk_1
\right)
dt', \quad \tfrac 12 < s < \tfrac 32.
\nn
}
Hence, noting that  
\ddd{
\int_{k_1\in\mathbb R} 
\no{R(k_1, t')}_{H^s(\mathbb R_{x_2})}^2
dk_1
\leqslant
\no{F(t')}_{H^s(\mathbb R_{x_1}\times \mathbb R_{x_2})}^2, \quad s\geqslant 0,
\nn
}
we infer the estimate
\ddd{ 
\no{W(x_2)}_{H_{x_1}^0 H_t^{\frac{2s+1}{4}}}
&\lesssim
\sup_{t\in [0, T]} \no{F(t)}_{H^s(\mathbb R_{x_1}\times \mathbb R_{x_2})}
\bigg[
\int_{t=0}^T \int_{z=0}^{T-t}
\frac{1}{z^{\frac 32 + s}}
\left(
\int_{t'=t}^{t+z}  dt'
\right)^2
dz dt
+
\int_{t'=0}^{T} dt'
\bigg]^{\frac 12}
\nn\\
&\lesssim
\sqrt{T} \sup_{t\in [0, T]} \no{F(t)}_{H^s(\mathbb R_{x_1}\times \mathbb R_{x_2})}, \quad \tfrac 12 < s < \tfrac 32, \ x_2\in\mathbb R.
\nn
}
In fact, starting from the bound \eqref{2d-nls-1d-fls-te-q-b} and proceeding as above, we obtain the same estimate also in the case $s=\frac 32$. Therefore, 
\eee{\label{2d-nls-fls-ivp-d1}
\no{W(x_2)}_{H_{x_1}^0 H_t^{\frac{2s+1}{4}}}
\lesssim
\sqrt{T} \sup_{t\in [0, T]} \no{F(t)}_{H^s(\mathbb R_{x_1}\times \mathbb R_{x_2})}, \quad \tfrac 12 < s \leqslant \tfrac 32, \ x_2\in\mathbb R.
}

Moreover, combining the bound \eqref{2d-nls-1d-fls-te-q-b} for $s=0$ with the norm \eqref{2d-nls-kat-1}, we find
\eee{
\no{W(x_2)}_{H_{x_1}^s H_t^{\frac 14}}^2
\lesssim
\int_{t=0}^T 
\int_{k_1\in\mathbb R}
\left(1+k_1^2\right)^s 
\no{e^{ik_1^2t} \what F^{x_1}(k_1, x_2, t)}_{L^2(\mathbb R_{x_2})}^2 dk_1 dt,
\nn
}
so by Parseval-Plancherel in $x_2$ and $k_2$ we infer
\ddd{\label{2d-nls-fls-ivp-d2}
\no{W(x_2)}_{H_{x_1}^s H_t^{\frac 14}}
&\lesssim
\left(\int_{t'=0}^T  \no{F(t')}_{H^s(\mathbb R_{x_1}\times \mathbb R_{x_2})}^2 dt'\right)^{\frac 12}
\nn\\
&\leqslant 
 \sqrt T \sup_{t\in [0, T]} \no{F(t)}_{H^s(\mathbb R_{x_1}\times \mathbb R_{x_2})}, \quad s\geqslant 0.
}

Estimates \eqref{2d-nls-fls-ivp-d1} and \eqref{2d-nls-fls-ivp-d2} imply estimate \eqref{2d-nls-fls-ivp-bst} via the definition  of the $B_T^s$-norm.
\end{Proof}

%
%
%
%
%
%
%
%
\section{The Forced Linear IBVP: Proof of Theorem \ref{2d-nls-fls-ibvp-t}}
\label{2d-nls-full-ibvp-s}

We shall now combine the pure IBVP Theorem \ref{2d-nls-ibvp-rv-t} with the IVP Theorems \ref{2d-nls-ls-ivp-t} and \ref{2d-nls-fls-ivp-t} in order to estimate the solution of the forced linear IBVP \eqref{2d-nls-fls-ibvp}. 
To do so, we first decompose this problem into simpler problems that have already been estimated in the preceding sections.

\vskip 3mm
\noindent
\textbf{Decomposition into simpler problems.}
IBVP \eqref{2d-nls-fls-ibvp} can be expressed as the superposition of  the homogeneous linear IBVP
\sss{\label{2d-nls-ls-ibvp}
\ddd{
&i u_t + u_{x_1x_1}+u_{x_2x_2}=0,
&&
(x_1, x_2, t)\in \mathbb R\times \mathbb R^+ \times (0, T),
\\
&u(x_1, x_2, 0) = u_0(x_1, x_2)\in H^s(\mathbb R_{x_1}\times \mathbb R_{x_2}^+), \quad && (x_1, x_2)\in \mathbb R\times \mathbb R^+, 
\\
&u(x_1, 0, t) = g_0(x_1, t) \in B_T^s, && (x_1,  t)\in \mathbb R\times [0, T],
}
}
and the forced linear IBVP with  zero  initial and boundary data
\sss{\label{2d-nls-fls-ibvp-0}
\ddd{
&i u_t + u_{x_1x_1}+u_{x_2x_2}= f \in C([0, T];  H^s(\mathbb R_{x_1}\times \mathbb R_{x_2}^+)),
\quad  &&
(x_1, x_2, t)\in \mathbb R\times \mathbb R^+ \times (0, T),
\\
&u(x_1, x_2, 0) = 0,  && (x_1, x_2)\in \mathbb R\times \mathbb R^+, 
\\
&u(x_1, 0, t) = 0, && (x_1,  t)\in \mathbb R\times [0, T].
}
}

The above decomposition has decoupled the forcing $f$ from the data $u_0, g_0$. Next, we will perform further decompositions in order to separate the data from each other. 

In particular, let the initial datum $U_0$ of the homogeneous linear IVP  \eqref{2d-nls-ls-ivp} be defined as a whole-plane extension of the half-plane initial datum $u_0$ of IBVP \eqref{2d-nls-ls-ibvp} such that
\eee{\label{2d-nls-uU}
\no{U_0}_{H^s(\mathbb R_{x_1}\times \mathbb R_{x_2})} \leqslant c\no{u_0}_{H^s(\mathbb R_{x_1}\times \mathbb R_{x_2}^+)}. 
}
Then, problem \eqref{2d-nls-ls-ibvp} can be expressed as the superposition of  IVP  \eqref{2d-nls-ls-ivp} and  the following homogeneous linear IBVP with  zero initial datum:
\sss{\label{2d-nls-ls-ibvp-r}
\ddd{
&i u_t + u_{x_1x_1}+u_{x_2x_2}=0,
&&
(x_1, x_2, t)\in \mathbb R\times \mathbb R^+ \times (0, T),
\\
&u(x_1, x_2, 0) = 0, && (x_1, x_2)\in \mathbb R\times \mathbb R^+, 
\\
&u(x_1, 0, t) = G_0(x_1, t) = g_0(x_1, t) - U(x_1, 0, t), \quad && (x_1,  t)\in \mathbb R\times [0, T],
}
}
where the function $U(x_1, 0, t)$ involved in the boundary condition is the solution $U=S\big[U_0; 0\big]$ of IVP \eqref{2d-nls-ls-ivp} evaluated at $x_2=0$ (this trace is well-defined since $s>\frac 12$).

In addition, let the forcing $F$ of the forced linear IVP  \eqref{2d-nls-fls-ivp} be defined as a whole-plane extension of the half-plane forcing $f$ of IBVP \eqref{2d-nls-fls-ibvp} such that
\eee{\label{2d-nls-fF}
\no{F}_{C([0, T]; H^s(\mathbb R_{x_1}\times\mathbb R_{x_2}))}
\leqslant
c  \no{f}_{C([0, T]; H^s(\mathbb R_{x_1}\times\mathbb R_{x_2}^+))}.
}
Then, problem \eqref{2d-nls-fls-ibvp-0} can be written as the superposition of   IVP  \eqref{2d-nls-fls-ivp} and the following homogeneous linear IBVP with  zero initial datum:
\sss{\label{2d-nls-fls-ibvp-r}
\ddd{
&i u_t + u_{x_1x_1}+u_{x_2x_2}=0,
&&
(x_1, x_2, t)\in \mathbb R\times \mathbb R^+ \times (0, T),
\\
&u(x_1, x_2, 0) = 0, && (x_1, x_2)\in \mathbb R\times \mathbb R^+, 
\\
&u(x_1, 0, t) = W(x_1, 0, t), \quad && (x_1,  t)\in \mathbb R\times [0, T],
}
}
where the boundary datum $W(x_1, 0, t)$ is obtained by evaluating the solution $W=S\big[0; F\big]$ of IVP \eqref{2d-nls-fls-ivp} at $x_2=0$ (this trace is well-defined since $s>\frac 12$). 

In summary, the solution of the forced linear IBVP \eqref{2d-nls-fls-ibvp} can be analyzed into the respective solutions of the four component problems \eqref{2d-nls-ls-ivp}, \eqref{2d-nls-fls-ivp}, \eqref{2d-nls-ls-ibvp-r} and \eqref{2d-nls-fls-ibvp-r} involved in the above decompositions as follows:
\eee{\label{2d-nls-fls-ibvp-utm-sol-T-sup}
S\big[u_0, g_0; f\big] 
=
S\big[U_0; 0\big]\Big|_{x_2 \in \mathbb R^+}
+
S\big[0; F\big]\Big|_{x_2 \in \mathbb R^+}
+
S\big[0, G_0; 0\big]
-
 S\big[0,  W\big|_{x_2=0}; 0\big].
}

The first two terms on the right-hand side of \eqref{2d-nls-fls-ibvp-utm-sol-T-sup} have been estimated in Theorems \ref{2d-nls-ls-ivp-t} and \ref{2d-nls-fls-ivp-t} respectively. 
Also, the remaining two terms can be estimated via Theorem \ref{2d-nls-ibvp-rv-t} as their associated problems, namely IBVPs \eqref{2d-nls-ls-ibvp-r} and \eqref{2d-nls-fls-ibvp-r},  essentially correspond to different versions of the pure IBVP \eqref{2d-nls-ibvp-rv}. Indeed, both of these problems are homogeneous, with zero initial datum, and with boundary datum in $B_T^s$   since for $\frac 12<s\leqslant \frac 32$ estimates \eqref{2d-nls-ls-ivp-te}, \eqref{2d-nls-fls-ivp-bst}, \eqref{2d-nls-uU} and \eqref{2d-nls-fF} imply
\ddd{\label{2d-nls-G0-est}
&\no{G_0}_{B_T^s}
\lesssim
\no{g_0}_{B_T^s}
+
\no{u_0}_{H^s(\mathbb R_{x_1}\times \mathbb R_{x_2}^+)},  
\\
\label{2d-nls-wb-est-h}
&\big\|W|_{x_2=0}\big\|_{B_T^s}
\lesssim
\sqrt T   \sup_{t\in [0, T]}   \no{f(t)}_{H^s(\mathbb R_{x_1}\times \mathbb R_{x_2}^+)}.
}
Moreover, thanks to the compatibility condition \eqref{2d-nls-comp-cond} and the initial conditions \eqref{2d-nls-lsivp-ic} and \eqref{2d-nls-flsivp-ic}, for $s>1$ the boundary data $G_0$ and $W|_{x_2=0}$ both vanish at $t=0$ for all $x_1\in\mathbb R$. 

Hence, restricting $1<s\leqslant \frac 32$, we can treat IBVPs \eqref{2d-nls-ls-ibvp-r} and \eqref{2d-nls-fls-ibvp-r}  simultaneously by considering the following problem:
\sss{\label{2d-nls-ibvp-r}
\ddd{
&iu_t+u_{x_1x_1}+u_{x_2x_2} = 0, &&(x_1, x_2, t)\in \mathbb R \times \mathbb R^+ \times (0,T),
\label{2d-nls-ibvp-r-eq} 
\\
&u(x_1, x_2, 0)= 0, && (x_1, x_2)\in \mathbb R \times \overline{\mathbb R^+},  
\label{2d-nls-ibvp-r-ic} 
\\
&u(x_1, 0, t) = Q_0(x_1, t)\in B_T^s,\quad &&(x_1, t)\in \mathbb R_{x_1} \times [0, T],
\label{2d-nls-ibvp-r-bc}
}
}
where the boundary datum satisfies the compatibility condition 
\eee{\label{2d-nls-h0-0}
Q_0(x_1, 0) = 0 \quad \forall x_1\in \mathbb R.
}

Next, we will reduce the estimation of IBVP \eqref{2d-nls-ibvp-r}  to that of the pure IBVP \eqref{2d-nls-ibvp-rv}. Let  
\eee{\label{2d-nls-q-def}
q_0(k_1, t) = e^{ik_1^2t} \, \what{Q_0}^{x_1}(k_1, t)
}
and note that $q_0\in H_t^{\frac{2s+1}{4}}(0, T)$ since $Q_0\in B_T^s$. Moreover, $q_0(k_1, 0) = 0$ for all $k_1\in\mathbb R$ thanks to condition \eqref{2d-nls-h0-0}. 
Let $E$ be an extension of $q_0$ from $[0, T]$ to $\mathbb R_t$ such that 
\eee{\label{2d-nls-qe}
\no{E(k_1)}_{H^{\frac{2s+1}{4}}(\mathbb R_t)}
\leqslant
c
\no{q_0(k_1)}_{H_t^{\frac{2s+1}{4}}(0,T)}, \quad c>0.
}
Consider the function 
$E_\theta(k_1, t) 
=
\theta(t) E(k_1, t) \in H^{\frac{2s+1}{4}}(\mathbb R_t)$, where  $\theta\in C^\infty_0(\mathbb R)$ with $\theta = 1$ on $[-1, 1]$, $\theta = 0$ on $(-2, 2)^c$ and $\no{\theta}_{L^\infty(\mathbb R)}=1$.
Since $T<1$,  we have $E_\theta  =q$ on $\mathbb R_{k_1}\times [0,T]$. Furthermore, combining the algebra property in $H^{\frac{2s+1}{4}}(\mathbb R_t)$ with estimate \eqref{2d-nls-qe}, we find
\eee{\label{2d-nls-gG-est}
\no{E_\theta(k_1)}_{H^{\frac{2s+1}{4}}(\mathbb R_t)} 
\leqslant  
c_s
\no{q_0(k_1)}_{H_t^{\frac{2s+1}{4}}(0, T)}, \quad k_1\in\mathbb R.
}
Also, $E(k_1, 0)=q_0(k_1, 0)=0$ and $\text{supp}(E_\theta)\subset \mathbb R_{k_1} \times (-2, 2)$ so that, in particular, $E_\theta(k_1, 2) = 0$.
Thus, $E_\theta \in H_0^{\frac{2s+1}{4}}(0, 2)$ and hence Theorem 11.4 of \cite{lm1972} implies that  the extension 
\eee{
\label{2d-nls-h-def-h}
q(k_1, t) 
= 
\left\{
\arraycolsep=4pt
\def\arraystretch{1}
\begin{array}{ll}
E_\theta(k_1, t), &t\in (0,2),
\\
0, &t\in (0,2)^c,
\end{array}
\right.
}
belongs to $H^{\frac{2s+1}{4}}(\mathbb R_t)$ with
$
\no{q(k_1)}_{H^{\frac{2s+1}{4}}(\mathbb R_t)} \leqslant c_s  \no{E_\theta(k_1)}_{H_0^{\frac{2s+1}{4}}(0,2)}
$ 
for all $k_1\in\mathbb R$ and $1<s \leqslant \tfrac 32$. 
Therefore, recalling estimate  \eqref{2d-nls-gG-est} we infer
\eee{\label{2d-nls-qq0}
\no{q(k_1)}_{H^{\frac{2s+1}{4}}(\mathbb R_t)} 
\leqslant
c_s
\no{q_0(k_1)}_{H_t^{\frac{2s+1}{4}}(0, T)},
\quad
1 < s \leqslant \tfrac 32, \ k_1\in\mathbb R,
}
Overall, we have constructed an extension $q$ of $q_0$ that satisfies the bound \eqref{2d-nls-qq0} and has support $\text{supp}(q)\subset \mathbb R_{k_1}\times (0, 2)$. 
Also, again by Theorem 11.4 of \cite{lm1972}, the extension $\widetilde q$ of $q$ by zero outside $[0, T]$ satisfies  
\eee{\label{2d-nls-bcdet-q-2}
\no{\widetilde q(k_1)}_{H^{\frac{1}{4}}(\mathbb R_t)}
\leqslant 
c
\no{q_0(k_1)}_{H_t^{\frac{1}{4}}(0, T)}, \quad k_1\in\mathbb R.
}
Thus, letting the boundary datum $g$ of the pure IBVP \eqref{2d-nls-ibvp-rv}  be given by
\eee{\label{2d-nls-gq}
\what{g}^{x_1}(k_1, t) = e^{ik_1^2t} q(k_1, t) = e^{ik_1^2t} \widetilde q(k_1, t), \quad (k_1, t)\in\mathbb R \times [0, 2],
}
and recalling \eqref{2d-nls-q-def}, we conclude that IBVP \eqref{2d-nls-ibvp-r} is embedded inside the pure IBVP \eqref{2d-nls-ibvp-rv}, that is 
$
S\big[0, Q_0; 0\big] = S\big[0, g; 0\big]\big|_{t\in [0, T]}.
$
Also, importantly, the above construction together with \eqref{2d-nls-q-def} and \eqref{2d-nls-gq} guarantees that
\eee{\label{2d-nls-g-bs}
\no{g}_{B^s} \lesssim \no{Q_0}_{B_T^s}, \quad 1<s\leqslant \tfrac 32.
}
Therefore, thanks to Theorem \ref{2d-nls-ibvp-rv-t} we deduce the estimate
\eee{\label{2d-nls-r-rv-est}
\sup_{t\in [0, T]} \no{S\big[0, Q_0; 0\big](t)}_{H^s(\mathbb R_{x_1}\times \mathbb R_{x_2}^+)} 
\lesssim
\no{Q_0}_{B_T^s}, \quad  1<s \leqslant \tfrac 32.
}

\vskip 2mm
\noindent
\textbf{Proof of Theorem \ref{2d-nls-fls-ibvp-t}.}
At this point, all four components of the superposition \eqref{2d-nls-fls-ibvp-utm-sol-T-sup} have been estimated. In particular, 
estimate \eqref{2d-nls-ls-ivp-se} of Theorem \ref{2d-nls-ls-ivp-t} and inequality \eqref{2d-nls-uU} imply 
\eee{
\sup_{t\in [0, T]} \no{S\big[U_0; 0\big](t)}_{H^s(\mathbb R_{x_1}\times \mathbb R_{x_2})}
\leqslant
c \no{u_0}_{H^s(\mathbb R_{x_1}\times \mathbb R_{x_2}^+)},  \quad  s\in\mathbb R,
\nn
}
estimate \eqref{2d-nls-fls-ivp-se} of Theorem \ref{2d-nls-fls-ivp-t} together with inequality \eqref{2d-nls-fF} and the fact that $T<1$ yield
\eee{
\sup_{t\in [0, T]} \no{S\big[0; F\big](t)}_{H^s(\mathbb R_{x_1}\times \mathbb R_{x_2})}
\leqslant
c \, \sqrt T \sup_{t\in[0, T]} \no{f(t)}_{H^s(\mathbb R_{x_1}\times \mathbb R_{x_2}^+)},  \quad s\in\mathbb R,
\nn
}
estimate \eqref{2d-nls-r-rv-est} and inequality \eqref{2d-nls-G0-est}  with $Q_0=G_0$  imply
\eee{
\sup_{t\in[0,T]} \no{S\big[0, G_0; 0\big](t)}_{H^s(\mathbb R_{x_1}\times \mathbb R_{x_2}^+)} 
\leqslant
c_s\, \Big(
\no{u_0}_{H^s(\mathbb R_{x_1}\times \mathbb R_{x_2}^+)} 
+ 
\no{g_0}_{B_T^s}
\Big), 
\quad 1 < s \leqslant \tfrac 32,
\nn
}
and estimate \eqref{2d-nls-r-rv-est} together with inequality \eqref{2d-nls-wb-est-h}  with $Q_0=W\big|_{x_2=0}$  yield
\eee{
\sup_{t\in[0,T]} \no{S\big[0, W\big|_{x_2=0}; 0\big](t)}_{H^s(\mathbb R_{x_1}\times \mathbb R_{x_2}^+)} 
\leqslant
c_s \sqrt T \sup_{t\in [0, T]} \no{f(t)}_{H^s(\mathbb R_{x_1}\times \mathbb R_{x_2}^+)}, \quad 1 < s \leqslant \tfrac 32.
\nn
}
Combining the above four estimates, we deduce estimate \eqref{2d-nls-fls-ibvp-se}  for the forced linear Schr\"odinger IBVP \eqref{2d-nls-fls-ibvp}.
\hfill $\blacksquare$

%
%
%
%
%
%
%
%
%

\section{Well-Posedness of NLS on the Half-Plane: Proof of Theorem \ref{2d-nls-nls-ibvp-wp-t}}
\label{2d-nls-lwp-s}

Using the forced linear IBVP estimate \eqref{2d-nls-fls-ibvp-se} and the contraction mapping theorem,  we shall now establish well-posedness of the NLS IBVP \eqref{2d-nls-nls-ibvp} in the sense of Hadamard, i.e. we shall show that there exists a unique solution to this problem that depends continuously on the initial and boundary data.

\vskip 3mm
\noindent
\textbf{Existence and uniqueness.}
Setting $f= \pm |u|^{p-1}u$ in the unified transform  formula \eqref{2d-nls-fls-ibvp-utm-sol-T} for the solution $S\big[u_0, g_0; f\big]$ of the forced linear IBVP \eqref{2d-nls-fls-ibvp} gives rise to the following iteration map for the NLS IBVP \eqref{2d-nls-nls-ibvp}:
\eee{\label{2d-nls-it-map}
u \mapsto \Phi u
=
\Phi_{u_0,g_0}u
=
S\big[u_0, g_0; \pm |u|^{p-1}u\big].
}
We shall show that this map  is a contraction in the space
\eee{\label{2d-nls-x-def}
X = C([0,T^*]; H^s(\mathbb R_{x_1}\times\mathbb R_{x_2}^+))
}
for some appropriate lifespan $T^*\in (0, T]$ to be determined.

\vskip 3mm
\noindent
\textit{\b{Showing that the map $u\mapsto \Phi u$ is onto $X$.}}
For $1<s\leqslant \frac 32$, estimate \eqref{2d-nls-fls-ibvp-se} implies
\eee{\nn
\no{\Phi u}_X
\leqslant
c_s
\Big(
\no{u_0}_{H^s(\mathbb R_{x_1}\times\mathbb R_{x_2}^+)}
+ 
\no{g_0}_{B_T^s}
+
\sqrt{T^*} \sup_{t\in[0,T^*]} \big\|\left(u\bar u\right)^{\frac{p-1}{2}}(t)\big\|_{H^s(\mathbb R_{x_1}\times\mathbb R_{x_2}^+)}
\Big),
}
where we have written $|u|^{p-1}=\left(u\bar u\right)^{\frac{p-1}{2}}$ since $\frac{p-1}{2}\in\mathbb N$. Thus, by repeated use of the algebra property in $H^s(\mathbb R_{x_1}\times\mathbb R_{x_2}^+)$ we obtain
\eee{\label{2d-nls-onto-est}
\no{\Phi u}_X
\leqslant 
c_s\, \Big(\no{u_0}_{H^s(\mathbb R_{x_1}\times\mathbb R_{x_2}^+)}+\no{g_0}_{B_T^s}+\sqrt{T^*} \no{u}_X^{p}\Big).
}
Let $B(0,r)=\left\{u\in X : \no{u}_X\leqslant r\right\}$ be a ball centered at $0$ with radius $r = 2c_s\no{(u_0,g_0)}_{D}= \no{u_0}_{H^s(\mathbb R_{x_1}\times\mathbb R_{x_2}^+)}+\no{g_0}_{B_T^s}$.
For $u\in B(0, r)$, estimate \eqref{2d-nls-onto-est} implies
\eee{
\no{\Phi u}_X
\leqslant 
\frac r2 + c_s  \sqrt{T^*} \, r^p.
\nn
}
Hence, choosing 
\eee{\label{2d-nls-lifespan-onto}
T^*
\leqslant
\frac{1}{4 c_s^2  r^{2(p-1)}}
}
ensures that $\Phi u\in B(0, r)$ whenever $u\in B(0, r)$.
\vskip 3mm
\noindent
\textit{\b{Showing that the map $u \mapsto \Phi u$ is  a contraction in $X$.}}
We shall show that   
\eee{\label{nls-hl-contraction-ineq}
\no{\Phi  u_1 - \Phi  u_2}_X\leqslant \frac 12 \left\| u_1-u_2\right\|_X
}
for any $u_1, u_2\in B(0,r)$ and $1 < s \leqslant \tfrac 32$. Noting that
\eee{\nn
\Phi  u_1 - \Phi  u_2
=
S\big[0, 0; \pm \left( |u_1|^{p-1}u_1 - |u_2|^{p-1}u_2 \right) \big]
}
and using estimate \eqref{2d-nls-fls-ibvp-se}, we obtain
\eee{\label{2d-nls-phiu1u2-0}
\no{\Phi  u_1 - \Phi  u_2}_X
\leqslant
c_s\,
 \sqrt{T^*} \sup_{t\in[0,T^*] }  \left\| \left(|u_1|^{p-1}u_1 -  |u_2|^{p-1}u_2\right)(t) \right\|_{H^s(\mathbb R_{x_1}\times\mathbb R_{x_2}^+)}.
}

Since $\tfrac{p-1}{2}\in\mathbb N$, we can write $|u_1|^{p-1}u_1 -  |u_2|^{p-1}u_2 =
\left(u_1 \bar u_1\right)^{\frac{p-1}{2}} u_1
- 
\left(u_2\bar u_2\right)^{\frac{p-1}{2}}  u_2,
$ 
which is a polynomial that vanishes when $u_1=u_2$. Hence, we must have
\eee{\label{2d-nls-nonlin-writing}
|u_1|^{p-1}u_1 -  |u_2|^{p-1}u_2
=
P(u_1, \bar u_1, u_2, \bar u_2) \left(u_1-u_2\right)
+
Q(u_1, \bar u_1, u_2, \bar u_2) \left(\overline{u_1-u_2}\right),
}
for some polynomials $P$ and $Q$ of degree $(p-1)$.

Combining inequality \eqref{2d-nls-phiu1u2-0} with the writing \eqref{2d-nls-nonlin-writing} and the algebra property in $H^s(\mathbb R_{x_1}\times\mathbb R_{x_2}^+)$, we find
\ddd{
\no{\Phi  u_1 - \Phi  u_2}_X
&\leqslant
c_s\,    \sqrt{T^*} 
\,
\Big(\no{P(u_1, \bar u_1, u_2, \bar u_2)}_X
+
\no{Q(u_1, \bar u_1, u_2, \bar u_2)}_X
\Big)
 \no{u_1-u_2}_X
\nn\\
&\leqslant
c_s \, \sqrt{T^*}  \, c_p \, r^{p-1}  
\left\| u_1 - u_2\right\|_X,  \quad c_p>1.
\nn
}
Hence, choosing
\eee{\label{2d-nls-lifespan-contr}
T^*
\leqslant
\frac{1}{4c_p^2c_s^2 r^{2(p-1)}}
}
ensures that the contraction inequality \eqref{nls-hl-contraction-ineq} is satisfied for all $u_1, u_2 \in B(0,r)$.

Note that, since $c_p >1$, the contraction condition \eqref{2d-nls-lifespan-contr} is stronger than the onto condition \eqref{2d-nls-lifespan-onto}. 
Overall, for $T^*\in (0, T]$ satisfying \eqref{2d-nls-lifespan-contr}, the map $u\mapsto\Phi u$ defined by  \eqref{2d-nls-it-map} is both onto and a contraction on the ball $B(0,r)$. Hence, by the contraction mapping theorem it follows that $u\mapsto\Phi u$ has a unique fixed point in $B(0,r)$. Equivalently, the integral equation $u=\Phi  u$  for the solution of the NLS IBVP \eqref{2d-nls-nls-ibvp} has a unique solution $u\in B(0,r) \subset X$. 

\vskip 3mm
\noindent
\textbf{Continuity of the data-to-solution map.} 
Finally, we shall show that the data-to-solution map
\eee{
H^s(\mathbb R_{x_1}\times\mathbb R_{x_2}^+)\times B_T^s
\ni
\left(u_0, g_0\right) 
\mapsto 
u\in X
}
is continuous.

Let $\left(u_0, g_0\right)$ and $\left(w_0, h_0\right)$ be two pairs of data  lying inside a ball $B_\varrho\subset D$ of radius $\varrho>0$ centered at a distance $r$ from $0$. 
Denote by $u=\Phi_{u_0, g_0} u$ and $w = \Phi_{w_0, h_0} w$  the  corresponding solutions to the NLS IBVP \eqref{2d-nls-nls-ibvp} and by $T_u$ and $T_w$ their respective lifespans, which are  given by 
\eee{\nn
T_u
=\min\left\{ T,  c_{s, p} \no{(u_0, g_0)}_{D}^{-2\left(p-1\right)}\right\},
\quad
T_w
=\min\left\{ T, c_{s, p} \no{(w_0, h_0)}_{D}^{-2\left(p-1\right)}\right\}, 
\quad
c_{s, p} =\left(2^{2p} c_p^2 c_s^{2p}\right)^{-1}. 
}
Since 
$
\max\left\{\no{\left(u_0, g_0\right)}_{D}, \no{\left(w_0, h_0\right)}_{D}\right\} \leqslant r+\varrho
$
and $p>1$,  it follows that
\eee{\label{nls-hl-Tc}
\min\left\{T_u, T_w\right\}
 \geqslant  
 \min\left\{ T,  c_{s, p} \left(r+\varrho\right)^{-2\left(p-1\right)} \right\}  = T_c.
 }
Hence, both solutions are guaranteed to exist for $0\leqslant t \leqslant T_c$. 

The common lifespan $T_c$ gives rise to the space
$
X_c = C([0, T_c]; H^s(\mathbb R_{x_1}\times\mathbb R_{x_2}^+)).
$
We shall now determine the radius $r_c$ of a ball $B(0, r_c)\subset X_c$ such that $u, w\in B(0, r_c)$ and  
\eee{
\label{nls-hl-lipaim}
\no{u - w}_{X_c} \leqslant 2c_s \left\| \left(u_0-w_0, g_0-h_0\right) \right\|_{D}.
}
Recall that $u$ and $w$ are  fixed points of the map $\Phi$ in the spaces $X_u$ and $X_w$  defined by \eqref{2d-nls-x-def} with $T^*$ replaced by $T_u$ and $T_w$ respectively. Hence, since $X_u, X_w\subset X_c$ by the definition of $T_c$,  $u-w$ is a fixed point of  $\Phi$ in $X_c$. Thus, using estimate \eqref{2d-nls-fls-ibvp-se}  together with the writing \eqref{2d-nls-nonlin-writing} and  the algebra property in $H^s(\mathbb R_{x_1}\times\mathbb R_{x_2}^+)$,  we obtain
\ddd{
\no{u - w}_{X_c}
&=
\left\| 
S\big[u_0-w_0,g_0-h_0, \pm \left(|u|^{p-1}u-w|w|^{p-1}\right) \big]
\right\|_{X_c}
\nn\\
&\leqslant 
c_s \left\| \left(u_0-w_0, g_0-h_0\right)\right\|_{D} 
+ 
c_s  c_p \sqrt{T_c}  \, r_c^{p-1}   \left\| u-w\right\|_{X_c}.
\label{2d-nls-ult}
}
For
$$
r_c
=
\big(2c_s c_p \sqrt{T_c}\,\big)^{-\frac{1}{p-1}}
= 
\max
\left\{
\big(2c_s c_p \sqrt{T}\,\big)^{-\frac{1}{p-1}},  
2c_s\left(r + \varrho \right)
\right\},
$$
estimate \eqref{2d-nls-ult} implies inequality  \eqref{nls-hl-lipaim} which in turn establishes local Lipschitz continuity of the data-to-solution map. The proof of Theorem \ref{2d-nls-nls-ibvp-wp-t} is complete.

\begin{remark}
Estimates \eqref{2d-nls-ibvp-rv-te}, \eqref{2d-nls-ls-ivp-te} and \eqref{2d-nls-fls-ivp-bst} give rise to the following additional estimate for the forced IBVP \eqref{2d-nls-fls-ibvp}:
$$
\sup_{x_2\in \overline{\mathbb R^+}} \no{u(x_2)}_{B_T^s} 
\leqslant 
c_s
\Big(
\no{u_0}_{H^s(\mathbb R_{x_1}\times\mathbb R_{x_2}^+)}
+ 
\no{g_0}_{B_T^s}
+
\sqrt T \sup_{t\in[0,T]} \no{f(t)}_{H^s(\mathbb R_{x_1}\times\mathbb R_{x_2}^+)}
\Big), \quad 1<s\leqslant \tfrac 32,
$$
which, as noted in Remark \ref{2d-nls-space-r}, can be combined with estimate \eqref{2d-nls-fls-ibvp-se} in order to show well-posedness of the NLS IBVP \eqref{2d-nls-nls-ibvp} in the smaller space $C([0, T^*]; H^s(\mathbb R_{x_1}\times \mathbb R_{x_2}^+))\cap C(\mathbb R_{x_2}^+; B_{T^*}^s)$ instead of the Hadamard space of Theorem \ref{2d-nls-nls-ibvp-wp-t}.
\end{remark}

%
%
%
%
%
%

\vskip 0.1in
\noindent
{\bf Acknowledgements.} This work was partially supported by a grant from the Simons Foundation (\#524469 to Alex Himonas).
The authors are grateful to Athanassios Fokas for being the catalyst of this project and for making valuable suggestions that led to the improvement of this work.

%
%
%
%
%
%
%
%
%

\vspace{8mm}
\noindent
A. Alexandrou Himonas \hfill Dionyssios Mantzavinos
\\
Department of Mathematics \hfill Department of Mathematics \\
University of Notre Dame   \hfill University of Kansas
\\
Notre Dame, IN 46556 \hfill Lawrence, KS 66045 \\
E-mail: \textit{himonas.1$@$nd.edu} \hfill E-mail: \textit{mantzavinos@ku.edu}

\vspace{3mm}

%
%
%
%
%
%
%
%
\section*{Appendix: Solution of the Forced Linear IBVP via Fokas' Unified Transform}
\renewcommand{\theequation}{A.\arabic{equation}}
\renewcommand{\thetheorem}{A.\arabic{theorem}}
\renewcommand{\thelemma}{A.\arabic{lemma}}
\renewcommand{\theremark}{A.\arabic{remark}}
\renewcommand{\thecorollary}{A.\arabic{corollary}}
\renewcommand{\thefigure}{A.\arabic{figure}}
\setcounter{corollary}{0}
\setcounter{equation}{0}
\setcounter{theorem}{0}
\setcounter{lemma}{0}
\setcounter{remark}{0}

We provide a concise derivation of the unified transform formula \eqref{2d-nls-fls-ibvp-utm-sol-T} for the forced linear Schr\"odinger IBVP \eqref{2d-nls-fls-ibvp} under the assumption of smooth initial and boundary values. For a more detailed derivation, we refer the reader to \cite{f2002}, where formula \eqref{2d-nls-fls-ibvp-utm-sol-T} was first obtained. 

Let $\widetilde u$ satisfy the adjoint of the linear Schr\"odinger equation, i.e.
$i\widetilde u_t - \widetilde u_{x_1x_1} - \widetilde u_{x_2x_2} = 0$.
Multiplying this equation by the solution $u$ of the forced linear Schr\"odinger equation \eqref{2d-nls-fls-eq}, and adding to it equation \eqref{2d-nls-fls-eq} multiplied by $\tilde u$,  we arrive at the divergence form
\eee{\nn
i\left(\widetilde u u\right)_t + \left(\widetilde u u_{x_1}-\widetilde u_{x_1}u\right)_{x_1} + \left(\widetilde u u_{x_2}-\widetilde u_{x_2}u\right)_{x_2}
=
\widetilde u  f. 
}
Setting $\widetilde u(x_1, x_2, t)=e^{-ik_1x_1-ik_2x_2 + i(k_1^2+k_2^2)t}$, $k_1, k_2\in\mathbb C$ and integrating with respect to  $x_1$ and $x_2$ while assuming that $u\to 0$ as $|x_1|, |x_2| \to \infty$, we obtain
\eee{\label{2d-nls-ode}
i\big(e^{i(k_1^2+k_2^2)t} \, \widehat u(k_1, k_2, t)\big)_t 
=
-e^{i(k_1^2+k_2^2)t}\left[\what{u}^{x_1}_{x_2}(k_1, 0, t) +ik_2 \what{u}^{x_1}(k_1, 0, t)\right] 
+ e^{i(k_1^2+k_2^2)t}  \widehat f(k_1, k_2, t),
}
where $\what u$ and $\what f$ are the half-plane Fourier transforms of $u$ and $f$ defined by \eqref{2d-nls-u0h-def} while $\what u^{x_1}$   denotes the Fourier transform  of $u$ with respect to $x_1$ defined as in \eqref{2d-nls-ft-x1-def}.

Note that  $\what u$ and $\what f$ are analytic as functions of $k_2$ in $\mathbb C^- = \left\{k_2\in\mathbb C:\text{Im}(k_2) < 0\right\}$ as a consequence of the following Paley-Wiener theorem.

\begin{theorem}[\b{\cite{s1994}, Theorem 7.2.4}]
\label{2d-nls-pw-t}
For any $\phi\in L^2(0, \infty)$, the half-line Fourier transform 
\eee{\nn
\what \phi(k) = \int_{x=0}^\infty  e^{-ikx} \phi(x) dx
}
is analytic in $\mathbb C^-= \left\{k\in\mathbb C: \textnormal{Im}(k)<0\right\}$.
\end{theorem}

Integrating \eqref{2d-nls-ode} with respect to $t$ gives rise to the so-called global relation
\ddd{\label{2d-nls-gr}
e^{i(k_1^2+k_2^2)t}\, \widehat u (k_1, k_2, t)
&=
\widehat u_0(k_1, k_2)
+
i\widetilde g_1(k_1, k_1^2+k_2^2, t) - k_2 \widetilde g_0(k_1, k_1^2+k_2^2, t)
\nn\\
&\quad
-i\int_{t'=0}^t e^{i(k_1^2+k_2^2)t'} \widehat f (k_1, k_2, t') dt',
\quad
(k_1, k_2)\in \mathbb R\times \overline{\mathbb C^-},
}
where for $g_j(x_1, t) = \p_{x_2}^j u(x_1, 0, t),$ $j=0, 1$, we define
\eee{\nn
\widetilde g_j(k_1, k_1^2+k_2^2, t) = \int_{t'=0}^t e^{i(k_1^2+k_2^2)t'} \what{g_j}^{x_1}(k_1, t')dt', 
\quad j=0,1.
}
Inverting \eqref{2d-nls-gr} by means of the usual inverse Fourier transform, we obtain  
\ddd{\label{2d-nls-ir}
&
u(x_1, x_2, t)
=
\frac{1}{(2\pi)^2} \int_{k_1\in\mathbb R}\int_{k_2\in\mathbb R} e^{ik_1x_1+ik_2x_2-i(k_1^2+k_2^2)t}\, \widehat u_0(k_1, k_2) dk_2 dk_1
\\
& 
-
\frac{i}{(2\pi)^2}\int_{k_1\in\mathbb R}\int_{k_2\in\mathbb R} e^{ik_1x_1+ik_2x_2-i(k_1^2+k_2^2)t} \int_{t'=0}^t e^{i(k_1^2+k_2^2) t'} \widehat f(k_1, k_2, t') dt' dk_2 dk_1
\nn\\
&
+
\frac{1}{(2\pi)^2} \int_{k_1\in\mathbb R}\int_{k_2\in\mathbb R} e^{ik_1x_1+ik_2x_2-i(k_1^2+k_2^2)t} \left[i\widetilde g_1(k_1, k_1^2+k_2^2, t) - k_2\widetilde g_0(k_1, k_1^2+k_2^2,t) \right]dk_2 dk_1.
\nn
}

\begin{figure}[ht]
\centering
\vspace{2.8cm}
\hspace{1cm}
\begin{tikzpicture}[scale=1.1]
\pgflowlevelsynccm
\draw[middlearrow={Stealth[scale=1.3]}, thick] (0,0) -- (90:2);
\draw[] (0,2) -- (0, 2.5);
\draw[dashed] (0,0) -- (0.5, 0);
\draw[dashed] (0,0) -- (0, -0.5);
\draw[middlearrow={Stealth[scale=1.3, reversed]},  thick] (0,0) -- (180:2);
\draw[] (-2,0) -- (-2.5,0);
\draw[middlearrow={Stealth[scale=1.1, reversed]}, thick, dashed] (-2, 0) arc (180:90:2);
\node[] at (-1.85,1.5) {\fontsize{10}{10} $\gamma_\rho^+$};
\node[] at (-2.3,-0.15) {\fontsize{8}{8} $-\rho$};
\node[] at (0.13, -0.23) {\fontsize{10}{10} $0$};
\node[] at (0.15,2.1) {\fontsize{8}{8} $i\rho$};
\end{tikzpicture}
\vspace{0.5cm}
\caption{Cauchy's theorem in the second quadrant of the complex $k_2$-plane.}
\label{2d-nls-deform-f}
\end{figure}
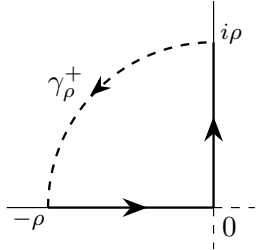

The representation \eqref{2d-nls-ir} is not an explicit solution formula for the forced linear IBVP \eqref{2d-nls-fls-ibvp} because it involves the unknown Neumann boundary value $u_x(0, t)$ through the transform $\widetilde g_1$. However, it turns out that $\widetilde g_1$ can be eliminated from \eqref{2d-nls-ir} in favor of known quantities.
In particular, note that since $x_2\geqslant 0$ and $t\geqslant t'$ the exponential $e^{ik_2x_2-i k_2^2(t-t')}$ is bounded for all $k_2 \in \overline{\mathbb C^+\setminus D}$, where $D$ here denotes the first quadrant of the complex $k_2$-plane. Thus, exploiting the analyticity of $\widetilde g_j$ for all $k_2\in\mathbb C$ (which follows via a Paley-Wiener theorem similar to Theorem \ref{2d-nls-pw-t}) we apply Cauchy's theorem in the second quadrant of the complex $k_2$-plane to write the $k_2$-integral in the last term of \eqref{2d-nls-ir} as
\eee{
\int_{k_2\in\p D} e^{ik_2x_2-ik_2^2t} \left[i\widetilde g_1(k_1, k_1^2+k_2^2, t) - k_2 \widetilde g_0(k_1, k_1^2+k_2^2, t) \right]dk_2
-
\lim_{\rho\to\infty} I(\rho, k_1, x_2, t),
\nn
}
where $\p D$ is the positively oriented boundary of  $D$ depicted in Figure \ref{2d-nls-dplus} and  for the quarter circle 
$
\gamma_\rho^+
=
\left\{\rho e^{i\theta}, \ \tfrac \pi 2\leqslant \theta \leqslant \pi\right\}
$
shown in Figure \ref{2d-nls-deform-f} we define
\eee{\label{2d-nls-Irho-def}
I(\rho, k_1, x_2, t)
=
 \int_{k_2\in \gamma_\rho^+} e^{ik_2x_2-ik_2^2t} \left[i\widetilde g_1(k_1, k_1^2+k_2^2, t) - k_2 \widetilde g_0(k_1, k_1^2+k_2^2, t) \right]dk_2.
}
Integrating by parts with respect to $t'$ in the definitions of $\widetilde  g_0$ and $\widetilde  g_1$, we obtain
\eee{\nn
\eqref{2d-nls-Irho-def}
\lesssim
\frac{C_{g_0, g_1}(k_1)}{x_2}
\cdot
\frac{ \rho \left(1-e^{-\rho x_2}\right)}{\left|\rho^2 - k_1^2\right|},
}
where, recalling that we work under the assumption of smooth boundary values,
\eee{\nn
C_{g_0, g_1}(k_1)
=
\no{\what{g_1}^{x_1}(k_1)}_{L_t^\infty(0, T)} + \no{\what{g_0}^{x_1}(k_1)}_{L_t^\infty(0, T)}
+
\no{\p_t \what{g_1}^{x_1}(k_1)}_{L_t^\infty(0, T)} + \no{\p_t\what{g_0}^{x_1}(k_1)}_{L_t^\infty(0, T)} <\infty.
}
Hence, the integral \eqref{2d-nls-Irho-def} vanishes in the limit $\rho\to \infty$ and, in turn, the integral representation  \eqref{2d-nls-ir} reads
\ddd{
\label{2d-nls-ir-def}
&u(x_1, x_2, t)
=
\frac{1}{(2\pi)^2} \int_{k_1\in\mathbb R}\int_{k_2\in\mathbb R} e^{ik_1x_1+ik_2x_2-i(k_1^2+k_2^2)t}\, \widehat u_0(k_1, k_2) dk_2 dk_1
\\
& 
-
\frac{i}{(2\pi)^2}\int_{k_1\in\mathbb R}\int_{k_2\in\mathbb R} e^{ik_1x_1+ik_2x_2-i(k_1^2+k_2^2)t} \int_{t'=0}^t e^{i(k_1^2+k_2^2) t'} \widehat f(k_1, k_2, t') dt' dk_2 dk_1
\nn\\
&
+
\frac{1}{(2\pi)^2} \int_{k_1\in\mathbb R}\int_{k_2\in\p D} e^{ik_1x_1+ik_2x_2-i(k_1^2+k_2^2)t} \left[i\widetilde g_1(k_1, k_1^2+k_2^2, t) - k_2\widetilde g_0(k_1, k_1^2+k_2^2,t) \right]dk_2 dk_1.
\nn
}

Next, note that under the transformation $k_2\mapsto -k_2$  the global relation \eqref{2d-nls-gr} yields the identity
\ddd{\label{2d-nls-gr-}
e^{i(k_1^2+k_2^2)t}\, \widehat u (k_1, -k_2, t)
&=
\widehat u_0(k_1, -k_2)
+
i\widetilde g_1(k_1, k_1^2+k_2^2, t) + k_2 \widetilde g_0(k_1, k_1^2+k_2^2, t)
\nn\\
&\quad
-i\int_{t'=0}^t e^{i(k_1^2+k_2^2)t'} \widehat f (k_1, -k_2, t') dt',
\quad
(k_1, k_2)\in \mathbb R\times \overline{\mathbb C^+}.
}
Thanks to the fact that
$$
\int_{k_2\in \p D}e^{ik_2x_2} \, \widehat u(k_1,-k_2,t)dk_2 = 0
$$
due to the analyticity and exponential decay of the integrand inside $D$, we are able to solve identity \eqref{2d-nls-gr-} for $\widetilde g_1$ and substitute the resulting expression in \eqref{2d-nls-ir-def} to obtain the explicit solution formula
\ddd{\label{2d-nls-fls-ibvp-utm-sol-t}
u(x_1, x_2, t)
&=
\frac{1}{(2\pi)^2}
\int_{k_1\in\mathbb R}
\int_{k_2\in\mathbb R}
e^{ik_1x_1+ik_2x_2-i(k_1^2+k_2^2)t}
\widehat u_0(k_1, k_2)
dk_2
dk_1
\nn\\
&
\quad
-\frac{1}{(2\pi)^2}
\int_{k_1\in\mathbb R}
\int_{k_2\in\p D}
e^{ik_1x_1+ik_2x_2-i(k_1^2+k_2^2)t}
\widehat u_0(k_1,-k_2)
dk_2
dk_1
\nn\\
&\quad
-\frac{i}{(2\pi)^2}
\int_{k_1\in\mathbb R}
\int_{k_2\in\mathbb R}
e^{ik_1x_1+ik_2x_2-i(k_1^2+k_2^2)t}
\int_{t'=0}^t 
e^{i(k_1^2+k_2^2)t'}\widehat f(k_1, k_2, t')dt'dk_2 dk_1
\nn\\
&
\quad
+\frac{i}{(2\pi)^2}
\int_{k_1\in\mathbb R}
\int_{k_2\in\p D}
e^{ik_1x_1+ik_2x_2-i(k_1^2+k_2^2)t}
\int_{t'=0}^t 
e^{i(k_1^2+k_2^2)t'}\widehat f(k_1, -k_2, t')dt'
dk_2 dk_1
\nn\\
&\quad
+\frac{1}{(2\pi)^2}
\int_{k_1\in\mathbb R}\int_{k_2\in\p D}
e^{ik_1x_1+ik_2x_2-i(k_1^2+k_2^2)t}
\, 2k_2 \, \widetilde g_0 (k_1, k_1^2+k_2^2, t)
dk_2dk_1.
}
Finally, exploiting once again analyticity and exponential decay in $D$, we infer that
$$
\int_{k_2\in \p D} e^{ik_2x_2-i k_2^2 t} \int_{t'=t}^T e^{i(k_1^2+k_2^2) t'} \what{g_0}^{x_1}(k_1, t') dt' dk_2 = 0
$$
and hence the solution formula \eqref{2d-nls-fls-ibvp-utm-sol-t} can be written in the equivalent form \eqref{2d-nls-fls-ibvp-utm-sol-T}, which is the one convenient for deriving linear estimates.

\end{document}